\documentclass[times]{nlaauth}

\usepackage{snapshot}

\newcommand{\A}{\mathrm{A}}
\newcommand{\C}{\mathrm{C}}
\newcommand{\E}{\mathrm{E}}
\newcommand{\I}{\mathrm{I}}

\newcommand{\M}{\mathrm{M}}

\newcommand{\V}{\mathrm{V}}
\newcommand{\W}{\mathrm{W}}
\newcommand{\X}{\mathrm{X}}
\newcommand{\Y}{\mathrm{Y}}

\newcommand{\mat}{\mathrm}
\newcommand{\bs}{\mathbf}
\renewcommand{\theta}{\vartheta}
\renewcommand{\phi}{\varphi}
\renewcommand{\rho}{\varrho}
\renewcommand{\epsilon}{\varepsilon}

\DeclareMathOperator{\Span}{span}
\DeclareMathOperator{\Range}{Im}

\DeclareMathOperator{\rank}{rank}
\usepackage[output-decimal-marker={,}]{siunitx}

 \usepackage[autostyle]{csquotes}

\usepackage[dvipsnames]{xcolor}

\usepackage{enumerate}
\usepackage[linesnumbered,ruled]{algorithm2e}
\usepackage{footnote}
\usepackage{hyperref}
\hypersetup{colorlinks=false}

\graphicspath{{./figures/}}
\DeclareGraphicsExtensions{.pdf}

\begin{document}

\title{Accelerating the iterative solution of convection-diffusion problems using singular value decomposition}

\runningheads{G.~Pitton and L.~Heltai}{Accelerating the iterative solution of convection-diffusion problems using SVD}

\author{G.~Pitton and L.~Heltai\corrauth}

\address{SISSA --- International School for Advanced Studies, via Bonomea 265, 34136 Trieste, Italy}

\corraddr{Luca Heltai \texttt{<luca.heltai@sissa.it>}}

\begin{abstract}

The discretization of convection-diffusion equations by implicit or semi-implicit methods leads to a sequence of linear systems usually solved by iterative linear solvers such as GMRES. Many techniques bearing the name of \emph{recycling Krylov space methods} have been proposed to speed up the convergence rate after restarting, usually based on the selection and retention of some Arnoldi vectors.

After providing a unified framework for the description of a broad class of recycling methods and preconditioners, we propose an alternative recycling strategy based on a singular value decomposition selection of previous solutions, and exploit this information in classical and new augmentation and deflation methods.
The numerical tests in scalar non-linear convection-diffusion problems are promising for high-order methods.

\end{abstract}

\keywords{GMRES; Augmentation Methods; Deflation Methods; Singular Value Decomposition}

\maketitle



\section{Introduction}
\label{sec:introduction}

The approximation of non-linear, time-dependent partial differential equations (PDE), independently on the actual numerical scheme used to discretise the equations and the domain, results in sequences of large and possibly non-symmetric linear systems, deriving from the linearisation of time varying non-linear terms. The most popular iterative solvers used for these problems are Krylov subspace methods, like the Conjugate Gradient (CG) method for symmetric positive definite systems~\cite{Hestenes1952}, the Generalized Conjugate Residual (GCR) method~\cite{Eisenstat1983}, and the Generalised Minimal Residual (GMRES) method~\cite{SaadSchultz} for general non-symmetric systems.

These are direct methods in exact arithmetic, since they find the exact solution of the system in at most $n$ iterations by computing an orthogonal basis for the entire vector space, and are approximate methods when the iteration process is terminated before convergence. In the CG method, the full Krylov subspace is explored through a three-term recursion, which guarantees a constant work and storage cost per iteration~\cite{GolubVanLoan}. This is in contrast with the GMRES method, where the  work and storage per iteration grow linearly with the iteration number. Despite the fact that finite precision arithmetic makes the construction of Krylov subspaces a more and more ill-conditioned process as the iteration count increases, it has been shown that the resulting loss of orthogonality among the Arnoldi vectors does not prevent the method from converging, as long as the vectors are linearly independent. The effect of finite precision on the stability of Krylov space methods was extensively studied in~\cite{GreenbaumRozloznik97,PaigeStrakos2002,PaigeRozloznik2006}.

In real life applications, the increasing storage cost of unrestarted GMRES makes it unfeasible for large non-symmetric systems. The most effective strategy to improve the rate of convergence in iterative solvers comes from the use of properly chosen preconditioners. While several preconditioners are available that perform very well for symmetric and positive definite linear systems, this is not the case for general non-symmetric systems, where the construction of good preconditioners is generally much more difficult and problem-dependent. When preconditioning alone fails to reduce the number of iterations below a reasonable threshold, it is common to use a restarted or truncated version of GMRES, or GMRES($m$). 

The most important side effect of the restart procedure is that the information of the vectors thrown away is lost at each restart, compromising or delaying the super-linear convergence that is often observed in the final GMRES iterations~\cite{VanderVorst1993}. The analysis of the GMRES(m) method~\cite{Joubert_On_1994} is non-trivial, and leads to some surprising results~\cite{Embree,GreenbaumPtak}. Some strategies vary heuristically the restart parameter in order to improve the convergence rate~\cite{Baker2009a}, but the most effective solutions try to keep orthogonality with respect to some previously spanned subspaces, for example by carefully retaining part of the previously constructed Krylov subspace, \emph{augmenting} the search space of the following cycle, or by \emph{deflating} the linear operator's spectrum~\cite{GCROT,SaadYeung2000}.

Even though these acceleration techniques are often presented as \emph{preconditioning} techniques, they are conceptually different. Deflation methods, first introduced in~\cite{Nicolaides} for the Conjugate Gradient method, consist in the left multiplication of the linear system by a suitably chosen projector, with the aim of removing the components responsible for slow convergence (e.g., small magnitude eigenvalues), resulting in a \emph{singular} algebraic system to be solved. In augmentation methods, a selected subspace responsible for a slowdown of the iterative method is removed from the matrix range, resulting in a ``tamed'' Krylov subspace.

We group these methods under the general name of \emph{acceleration methods}, as they aim at accelerating the computation of an approximate solution by using wisely at least two vector spaces: a Krylov space and some recycling space. The first one is usually obtained through Arnoldi or two-sided Lanczos algorithms, while the generation and the exploitation of the second vector space is largely algorithm-dependent. An application of these ideas to sequences of linear systems is presented in~\cite{deSturler2006}, and for Computational Fluid Dynamics (CFD) applications in~\cite{Amritkar_Recycling_2015}. 
Recycling methods for sequences of symmetric linear systems arising in topology optimization are studied in~\cite{Wang2007}, and in~\cite{Mello2010} and~\cite{Kilmer2006} are applied as efficient solvers for inverse problems in electrical impedance tomography.

However, when solving time-dependent problems, one may expect that it should be possible to extract useful spectral information also from a number of solutions to previous time steps, and not only from previous GMRES cycles, even though the solutions are obtained from different matrices and right-hand sides. 

We provide a unified framework for six different recycling methods, some of which are equivalent to deflation preconditioners known in the literature. These recycling methods are applied on the basis of subspaces that exploit the information contained in a number of previous solution steps (using Singular Value Decomposition, as done first in~\cite{Ali2012}) to accelerate the solution of the current time step.

In particular, we explore recycling methods in two different variations, by i) forming good initial guesses for the solution of the current linear systems, and by ii) constructing subspaces for augmented and deflated methods. The two approaches are known to be related to each other, and are often interpreted similarly.

The first procedure is inspired by extrapolation methods, commonly used in time-dependent problems to extrapolate the solution at the current time  to compute non-linear terms of the system. Extrapolated solutions can be used as a seed for Krylov subspaces, and in this work we show that a much more efficient result is obtained by replacing extrapolation with projection, as in~\cite{Ali2012}, and combining it with proper recycling Krylov methods. This removes the stability issues associated to equi-spaced extrapolation, allowing one to keep a larger number of previous solutions, and to obtain best approximation-type estimates.

Moreover, by replacing extrapolation with projection, additional freedom is available in the choice of the projection subspace. In this work we explore and then compare several existing possibilities, starting from orthonormalization of some previous solutions, and propose a new approach based on the retainment of the principal components of a set of previous solutions using Singular Value Decomposition (SVD).

A similar strategy is used to construct efficient augmented and deflative methods, that exploit the SVD of a collection of previous time step solutions to construct the recycling space. We use this space to construct several preconditioners and show that in many test cases, applying SVD for the selection of the recycling spaces works very well, at least for scalar convection-diffusion problems, and for moderate values of the P\'{e}clet number.

In Sections~\ref{sec:petrov-galerkin} and~\ref{sec:kryl-subsp} we introduce the setting of our time-dependent problem in terms of general Petrov--Galerkin approximations, and review the solution of time-dependent linear problems in terms of Krylov subspace methods, along with some of the most popular acceleration methods. We introduce a new framework for the unified discussion of some recycling spaces, that allows the systematic construction of three augmentation and three deflation methods, some of which are equivalent to known methods, and some of which are new. In Section~\ref{sec:accel-init-guess} we consider procedures for the construction of a good initial guess based on extrapolation or restriction on a space of previous solutions, and we propose a new technique based on a Multiple Subspace Correction principle, while in Section~\ref{sec:construction_augmentation} we review some techniques for the generation of a recycling subspace, and we consider two cases where the subspace is obtained from a fixed set of previous solutions, or by means of an SVD performed on the space spanned by some previous solution vectors. In Sections~\ref{sec:numerical_experiments} and~\ref{sec:conclusions} we present a selection of numerical experiments and draw some conclusions. 

\section{Problem definition}
\label{sec:petrov-galerkin}

Let $X$ and $Y$ be normed spaces, $Y'$ the dual of $Y$, and $L(u):X\to Y'$ a time-dependent, linear, and bounded map, possibly depending on the current solution $u\in X$, that represents the differential operator of a partial differential equation (PDE) on the space $X$. We are interested in finding approximate solutions of non-linear time-dependent PDEs of the form
\begin{equation}
  \partial_tu + L(u) u = f,
  \label{eq:pde}
\end{equation}
where $\partial_tu$ is a shorthand for the time derivative of $u$, and $f\in Y'$ is given.

In the following, we will consider principally a weak form of equation~\eqref{eq:pde}:
\begin{equation}
  \langle v,\partial_tu\rangle + \langle v,L(u) u\rangle = \langle v,f\rangle\qquad \forall v\in Y.
  \label{eq:weak_pde}
\end{equation}

From the weak equation~\eqref{eq:weak_pde} we can derive a Petrov--Galerkin approximation~\cite{Quarteroni} by introducing some finite dimensional subspaces $X_h\subset X$ and $Y_h\subset Y$ of dimension $n$, so that the approximation $u_h\in X_h$ solves:
\begin{equation}
  \langle v,\partial_tu_h\rangle + \langle v,L(u_h) u_h\rangle = \langle v,f\rangle\qquad \forall v\in Y_h.
  \label{eq:petrov-galerkin}
\end{equation}

It is understood that, as the discretization parameter $h$ tends to $0$, $X_h$ tends to span all of $X$. Letting $\{\phi_j\}_{j=1}^n$ and $\{\psi_i\}_{i=1}^n$ be two bases sets for $X_h$ and $Y_h$ respectively, the approximate solution can be written as: $u_h=\sum_{j=1}^n u_j(t)\phi_j$. Applying a finite-difference like discretization of the time derivative, the coefficients $u_j(t)$ are evaluated only at certain fixed points in time $u_j(t_k)$. We obtain a sequence of non-linear problems:
\begin{multline}
  \beta_0\langle\psi_i,\phi_j\rangle u_j(t_k) + \sum_{j=1}^n\langle \psi_i,L(u_h)\phi_j\rangle u_j(t_k) = \\
  \langle \psi_i,f\rangle - \sum_{\alpha=1}^p\beta_\alpha\langle\psi_i,\phi_j\rangle u_j(t_{k-\alpha})\qquad i=1,\dots,N,
  \label{eq:linear_pde}
\end{multline}
where $\beta_\alpha$ are the coefficients coming from the chosen finite difference scheme and $p$ is the number of steps of the time advancing scheme.

The non-linear term can be linearised around a (known)  guess solution $\tilde u$ by considering, for example, a Newton-like approximation 
\begin{equation}
 L(u_h)\sim L(\tilde u_h) + D_uL(\tilde u_h) (u_h-\tilde u_h),
\end{equation}
where $D_u L$ is the Fr\'{e}chet derivative of $L$. 

This leads to the time-dependent linear systems:
\begin{equation}
  \A(t_k) \bs{x}=\bs{b}(t_k),
  \label{eq:algebra_pde}
\end{equation}
where
\begin{equation}
  \begin{aligned}
(\A(t_k))_{ij} &= \beta_0\langle\psi_i,\phi_j\rangle + \langle\psi_i,L(\tilde u_h)\phi_j\rangle + \langle\psi_i,D_u L(\tilde u_h)\phi_j\rangle\\
(\bs{x})_i &= u_i(t_k)\\
(\bs{b})_i&=\langle\psi_i,f\rangle+\sum_{\alpha=1}^p\beta_\alpha\langle\psi_i,\phi_j\rangle u_j(t_{k-\alpha})+\langle\psi_i,D_u L(\tilde u_h) \tilde u_h \rangle.
\end{aligned}
  \label{eq:def_Axb}
\end{equation}

In the examples presented in this work, we neglect the term $ D_uL(\tilde u_h) (u_h-\tilde u_h)$, and consider simple extrapolations for $\tilde u_h$.  For the sake of simplicity, we assume that $\A\in\mathbb{R}^{n\times n}$, $\bs{x},\bs{b}\in\mathbb{R}^n$, but in principle it is possible to recast the problem for complex-valued functions by replacing the transpose operation with the hermitian and by using duality pairings in normed spaces over $\mathbb{C}$.

In the rest of the current section, we consider mainly the linear system:
\begin{equation}
\A\mathbf{x}=\mathbf{b}
\label{eq:sequence_linear_systems}
\end{equation}
where, to simplify the notation, the time-dependence of $\A$ is left implicit.
In general, the dimension $n$ of the linear system can depend on $t_k$, but this is not emphasized here to avoid the proliferation of indices. In general, we assume that the sequence of matrices and right-hand sides depend on the previous solution, $\A(t_i)=\A(t_i, \mathbf{x}(t_{i-1})),\bs{b}^i=\bs{b}^i(t_i, \bs{x}(t_{i-1}))$, and they are not entirely available until the very last element of the sequence.

\section{Review of Krylov subspace methods}
\label{sec:kryl-subsp}

Given an initial guess $\mathbf{x}_0$ for each timestep, Krylov subspace methods for linear systems (see, for instance,~\cite{GolubVanLoan,Greenbaum,Hackbusch,LiesenStrakos,Saad}) provide a sequence of approximate solutions $\{\mathbf{x}_k\}$ in the form (following the notation of~\cite{LiesenStrakos}):
\begin{equation}
 \begin{split}
   \mathbf{x}_k &= \mathbf{x}_0 + \bs{z}_k \\
   \mathbf{r}_k & = \mathbf{r}_0 - A\bs{z}_k\qquad
  \end{split}
  \qquad \bs{z}_k\in\mathcal{S}_k,
\end{equation}
where $\mathcal{S}_k$ is a $k-$dimensional \emph{search space}, and the approximation vector $\bs{z}_k$ is usually determined by imposing some constraints on the residual $\mathbf{r}_k:=\mathbf{b}-\A\mathbf{x}_k$.

Two common choices to select $\bs{z}_k$ are:
\begin{description}
\item[minimal residual]: $\min_{\mathbf{z}\in\mathcal{S}_k}\|\bs{r}_k\|$;
\item[orthogonal residual]: $\mathbf{r}_k\perp \mathcal{C}_k$ for some \emph{constraint space} $\mathcal{C}_k$.
\end{description}
  
We denote with $\{\bs{c}_i\}_{i=1}^k$ and $\{\bs{s}_i\}_{i=1}^k$ two basis sets for $\mathcal{C}_k$ and $\mathcal{S}_k$, that will be identified with the column vectors of the orthogonal matrices $\mat{C}_k\in\mathbb{R}^{n\times k}$ and $\mat{S}_k\in\mathbb{R}^{n\times k}$. 

In practice, if $\mat{S}_k,\C_k\in\mathbb{R}^{n\times k}$ are constructed such that their column space spans respectively the search space and the constraint space, i.e. $\Range{\mat{S}_k}=\mathcal{S}_k$ and $\Range{\C_k}=\mathcal{C}_k$, then the orthogonal residual solutions can be characterized by:
\begin{equation}
 \C_k^T\mathbf{r}_k = 0 \quad \Longrightarrow  \quad \C_k^T\A\mat{S}_k\mathbf{t}_k=\C_k^T\mathbf{r}_0,
\label{eq:reduced_linear_system}
\end{equation}
where $\mathbf{z}_k = \mat{S}_k\mathbf{t}_k$, $\mathbf{t}_k\in\mathbb{R}^k$, and $\mathbf{r}_0=\mathbf{b}-\A\mathbf{x}_0$ is the initial residual.
In this case, the approximate solution at the $k$-th iteration is given by:
\begin{equation}
  \begin{split}
    \mathbf{x}_k &=\mathbf{x}_0+\mat{S}_k(\C_k^T\A\mat{S}_k)^{-1}\C_k^T\mathbf{r}_0\\
    \mathbf{r}_k &=\mathbf{r}_0-\mat{A}\mat{S}_k(\C_k^T\A\mat{S}_k)^{-1}\C_k^T\mathbf{r}_0.
  \end{split}
\label{eq:x_k_expression}
\end{equation}

While minimal residual solutions exist and are unique for any choice of $\mathcal{S}_k$, this is not the case for the orthogonal residual case. Indeed, the existence of orthogonal residual solutions can be ensured provided that some compatibility conditions hold between the search and constraint spaces. However,  for some choices of $\mathcal{S}_k$, the minimal residual and orthogonal residual solutions coincide. We refer to~\cite{LiesenStrakos} for an in-depth discussion on this topic.

Equation~\eqref{eq:x_k_expression} can be expressed in terms of a projection operator:
\begin{equation}
  \begin{split}
    \bs{x}_k &=\bs{x}_0+\A^{-1}\mat{P}_{\A\mathcal{S}_k,\mathcal{C}_k}\bs{r}_0,\\
    \bs{r}_k &=(\mat{I}- \mat{P}_{\A\mathcal{S}_k,\mathcal{C}_k})\bs{r}_0,
  \end{split}
\label{eq:def_x_k_A_P}
\end{equation}
where the notation $\mat P_{\mathcal{A},\mathcal{B}} \bs x$ denotes the projection of $\bs x$ onto the space $\mathcal{A}$ orthogonal to $\mathcal{B}$, i.e.,  $\mat{P}_{\A\mathcal{S}_k,\mathcal{C}_k} :=\A\mat{S}_k(\C_k^T\A\mat{S}_k)^{-1}\C_k^T$ is a projector onto $\A\mathcal{S}_k$ orthogonal to $\mathcal{C}_k$ (see e.g.~\cite{LiesenStrakos}). 

We remark that in all practical applications, the iteration $\mathbf{x}_k$ is not computed naively by solving the reduced linear system~\eqref{eq:reduced_linear_system} but in more computationally efficient ways~\cite{Saad}.

Now we are left with the choice of the spaces $\mathcal{S}_k$ and $\mathcal{C}_k$. The most popular methods construct their iterations on Krylov subspaces, i.e., successive powers of $\mat{A}$ times the initial residual:
\begin{equation}
\mathcal{S}_k = \mathcal{K}_k(\A,\mathbf{r}_0):=\Span\{\mathbf{r}_0,\A\mathbf{r}_0,\A^2\mathbf{r}_0,\dots,\A^{k-1}\mathbf{r}_0\},
\end{equation}
and in the nonsymmetric case the most frequent choices are $\mathcal{C}_k=\mathcal{S}_k$ for the FOM, $\mathcal{C}_k=\A\mathcal{S}_k$ for the GMRES method~\cite{SaadSchultz} and $\mathcal{C}_k=\mathcal{K}_k(\A^T,\mathbf{r}_0)$ for the methods based on the nonsymmetric Lanczos procedure such as BiCGStab and the QMR family.

In Krylov space methods, two directions can be explored to speedup the convergence rate, namely:
\begin{enumerate}[i)]
  \item select a good initial guess, so that some norm of $\mathbf{r}_0$ is as small as possible;
  \item manipulate the matrix (and possibly the right-hand side) in order to obtain an equivalent linear system that is easier to solve.
\end{enumerate}
These two strategies are not always so sharply distinguished, as is the case with multiple subspace correction methods.

Manipulation of the system matrix is usually referred to as \emph{preconditioning} of the system (see \cite{Benzi} for a thorough introduction): both sides of the equation could be multiplied on the left by a matrix $\M$ (called in this case \emph{left preconditioner}) such that the new linear system:
\begin{equation}
\M\A\mathbf{x}=\M\mathbf{b}
\end{equation}
can be solved more effectively by the iterative method. 
An alternative is \emph{right preconditioning}:
\begin{equation}
\A\M\mathbf{y}=\mathbf{b}\qquad \mathbf{x}=\M^{-1}\mathbf{y},
\end{equation}
that in some circumstances allows much more freedom in the choice of $\M$~\cite{SaadFlexible}.

When both approaches are used at the same  time, it is desirable to make sure that the resulting Krylov subspaces do not contain components in the directions that have already been explored when constructing the initial guess, resulting in multi-level preconditioners.

Over the last 20 years, many techniques have been proposed to exploit this principle in restarted GMRES, most of them under the name of \emph{recycling Krylov spaces}, \emph{augmentation methods} or \emph{deflation methods}. A common feature shared by these \emph{acceleration methods} is to compute approximate solutions to large linear systems using wisely at least two low-dimensional vector spaces: a Krylov space and some recycling space. The first one is usually obtained through Arnoldi or two-sided Lanczos algorithms, while the generation and the exploitation of the second vector space is largely algorithm-dependent.

\subsection{Augmented Krylov methods}
\label{sec:augmented_krylov}

The discussion starting in this section aims at providing a common framework for some augmentation and deflation methods, described in Section~\ref{sec:synoption_aug_def}. The theoretical treatment of recycling methods is kept to a minimum, referring to~\cite{Saad97} and~\cite{ChapmanSaad97} for a general analysis of deflated and augmented methods.

Augmentation methods enrich the search space $\mathcal{S}_k$ of the Krylov method with an ``augmentation subspace'' $\mathcal{V}_s$, so that the approximate solution will lie in the new search space:
\begin{equation}
  \mathcal{S}_k=\mathcal{K}_k(\A,\mathbf{r}_0)+\mathcal{V}_s.
  \label{eq:def_aug_S_space}
\end{equation}
As a consequence, in augmentation methods applied to the GMRES solver, the residual at the $k$-th step is orthogonal to the following ``augmented'' constraint space:
\begin{equation}
  \mathcal{C}_k=\mathcal{K}_k(\A,\A\mathbf{r}_0)+\A\mathcal{V}_s.
  \label{eq:def_aug_C_space}
\end{equation}
This method is effective in reducing the iteration count when the space $\mathcal{V}_s$ contains a good approximation of the solution, and when successive  powers $A^k\mathbf{r}_0$ fail to cover the elements of $\mathcal{V}_s$ within a reasonable number of iterations. Due to the additional orthogonalization operations required by larger search and constraint spaces, the dimension $s$ of $\mathcal{V}_s$ should be chosen with a grain of salt to achieve a good compromise in the computational cost, as will be explained in Section~\ref{sec:efficiency}.

In the case of augmentation methods, the $k$-the iteration does not formally differ from the one in equation~\eqref{eq:x_k_expression}, the only difference being the specific choices~\eqref{eq:def_aug_S_space} and~\eqref{eq:def_aug_C_space} for the matrices $S_k$ and $C_k$.

\subsection{Deflated Krylov methods}
\label{sec:adapted_deflation}

If the recycling  subspace $\mathcal{V}_s$ does not contain a good approximation of the solution, i.e., the solution projection onto this space has a small magnitude, but at the same time  $\mathcal{V}_s$  has components along eigenvectors of $\A$ with small eigenvalues, it may be convenient to \emph{deflate} the search space of the entire subspace  $\mathcal{V}_s$ (after solving for the components of the solution along $\mathcal{V}_s$).


In the following,  we denote by $\V_s$ any orthogonal matrix whose columns span the space $\mathcal{V}_s$, and by $\E_s := \V_s^T \A \V_s$ the restriction of $\A$ onto $\mathcal{V}_s$, and we consider only left preconditioning (the right preconditioning counterparts are constructed similarly).

The linear system is multiplied on the left by the projector $\M_D=\I-\A\V_s\E_s^{-1}\V_s^T$, and the resulting problem:
\begin{equation}
  \M_D\A\bs{x}=\M_D\bs{b}.
  \label{eq:left_deflation}
\end{equation}
is solved with a Krylov space method.

The ``preconditioner'' qualifier for the deflation projection $\M_D$ should be taken with care, since the resulting operator is no longer invertible, and many authors make a distinction in this case, avoiding the term ``preconditioner''. The effect of $\M_D$ on $\A$ is to set to zero any element of $\mathcal{V}_s$, i.e., by construction $\M_D \A \V_s = 0$.


More generally, deflation preconditioners can be used to move some eigenvalues of the matrix $A$ to any target value $\lambda^*$, or to inject artificial eigenvalues into the resulting system by adding a term to the preconditioner:
\[
  \M_D=\I-\A\V_s\E_s^{-1}\V_s^T + \lambda^*\V_s\E_s^{-1}\V_s^T.
\]

Standard deflation preconditioners construct $\V_s$ using (an approximation of) the eigenvectors corresponding to the smallest eigenvalues of $\A$.

The references~\cite{NabbenVuik} and~\cite{Tang} present a comparison of many deflation preconditioner instances with standard multigrid and domain decomposition preconditioners. In~\cite{NabbenVuik}, the authors prove that the deflative preconditioners perform better than coarse grid correction preconditioners for any choice of the deflation space $\V_s$.
This is contrasted by the class of Adaptive Algebraic Multigrid~\cite{FalgoutAAM} and Adaptive Smoothed Aggregation Multigrid methods~\cite{Falgout_aSA,Falgout_aSAM}, where a coarse grid correction method in the form of an Algebraic Multigrid solver/preconditioner is enriched with the information obtained from relevant algebraic subspaces to enhance the robustness and to extend the range of applicability of Algebraic Multigrid methods.
The founding idea of Adaptive Smoothed Aggregation Multigrid is that of retrieving information on the subspace where the smoother is less effective by solving the homogeneous equation associated to the original problem.
The solution of the homogeneous equation provides an approximation for the ``coarse'' components of the error, which are then addressed by selecting the coarser correction spaces precisely with the goal of reducing those error components on which the smoother has proven to be not very effective.
An additional advantage of Adaptive Multigrid methods is the possibility of incorporating some known information coming e.g. from some recycling strategy in the definition of the hierarchy of coarse spaces and of the related prolongation/restriction operators.
This feature could be the starting point for an application of Adaptive Multigrid methods to the solution of sequences of linear systems, coupling these methods with some recycling technique to retain the relevant information coming from the solution of the previous linear systems.
Such possibility however is conceptually quite different from the approaches used in standard augmentation/deflation preconditioners for Krylov Space Methods, and is left for future investigation.

Finally, in~\cite{Tang} it is shown that for symmetric positive definite systems, deflation, abstract balancing and two-level multigrid preconditioners are comparable in terms of work per iteration and convergence properties. In the following we consider only standard deflation preconditioners.

We rewrite the $k$-th iteration as:
\begin{equation}
  \bs{x}_k=\bs{x}_0 + \mat{S}_k(\C_k^T(\I-\A\V_s\E_s^{-1}\V_s^T)\A\mat{S}_k)^{-1}\C_k^T(\I-\A\V_s\E_s^{-1}\V_s^T)\bs{r}_0.
  \label{eq:left_deflation_iterate}
\end{equation}
If the initial guess is built as for the multiple subspace correction case, namely $\bs{x}_0=\V_s\E_s^{-1}\V_s^T\bs{b}$, equation~\eqref{eq:left_deflation_iterate} becomes:
\begin{equation}
  \bs{x}_k=\V_s\E_s^{-1}\V_s^T\bs{b} + \mat{S}_k(\C_k^T(\I-\A\V_s\E_s^{-1}\V_s^T)\A\mat{S}_k)^{-1}\C_k^T(\I-\A\V_s\E_s^{-1}\V_s^T)\bs{b},
  \label{eq:left_deflation_special}
\end{equation}
where it can be immediately seen that the difference between the deflative approach and the Multiple Subspace Correction Method consists in the form of the projector acting on the left of the system matrix.
In this case, we see that the solution is again obtained through a projection of the initial residual, except that now it cannot be seen as an instance of the Multiple Subspace Correction Method. Indeed, equation~\eqref{eq:left_deflation_iterate} clearly shows the projector $\mat{P}_{\mathcal{V}_s,\A\mathcal{V}_s}=\A\V_s\E_s\V_s^T$ onto $\mathcal{V}_s$ orthogonal to $\A\mathcal{V}_s$, but this time the projection onto $\mathcal{S}_k$ is modified so that in the new operator:
\[
  \mat{P}_{\A\mathcal{S}_{k},\mathcal{C}_k,D}=\A\mat{S}_k(\C_k^T(\I-\A\V_s\E_s^{-1}\V_s^T)\A\mat{S}_k)^{-1}\C_k^T
\]
the deflation space is explicitly removed from the range of $A$ by means of the projector $(\I-\mat{P}_{\mathcal{V}_s,\A\mathcal{V}_s})\A$, so that the ``troublesome'' eigenvalues are removed from the spectrum of the modified operator.

Regarding the convergence properties of the deflation preconditioner compared to the augmentation strategies, in~\cite{EiermannErnstSchneider} it is shown that for minimal residual approximations, if the augmentation/deflation space $\mathcal{V}_s$ spans an invariant subspace\footnote{meaning that $\A\mathcal{V}_s\subseteq\mathcal{V}_s$} of $\A$, then
\begin{equation}
  0=\|\mat{P}_{\mathcal{V}_s}\bs{r}_{\mathrm{aug}}\| \le \|\mat{P}_{\mathcal{V}_s}\bs{r}_{\mathrm{def}}\| \qquad\text{and} \qquad \|(\I-\mat{P}_{\mathcal{V}_s})\bs{r}_{\mathrm{aug}}\| \le \|(\I-\mat{P}_{\mathcal{V}_s})\bs{r}_{\mathrm{def}}\|.
  \label{eq:estimate_deflation_vs_augmentation}
\end{equation}

\subsection{A unified view of augmented and deflated methods}
\label{sec:synoption_aug_def}

Augmented Krylov methods described in section~\ref{sec:augmented_krylov} and deflation preconditioners of section~\ref{sec:adapted_deflation} can be interpreted as Krylov space methods with a special choice of the initial guess, followed by suitable modifications of the Krylov space, usually achieved by projections.
We suppose that given a recycling space $\V_s$, the standard initial guess is obtained as in equation~\eqref{eq:projection_t_s}, namely $\bs{x}_0=\V_s\E_s^{-1}\V_s^T\bs{b}$. Writing the linear system's solution as the sum of the initial guess $\bs{x}_0$ and a correction $\bs{x}_c$, $\bs{x}=\bs{x}_0+\bs{x}_c$, we get the equation for the correction:
\begin{equation}
  \A\bs{x}_c=\bs{r}_0=(\I-\mat{P}_{\mathcal{V}_s,\A\mathcal{V}_s})\bs{b}.
  \label{eq:correction_equation}
\end{equation}
Solving equation~\eqref{eq:correction_equation} with a Krylov method leads to the search space
\begin{equation}
  \mathcal{S}_k=\mathcal{K}_k(\A,(\I-\mat{P}_{\mathcal{V}_s,\A\mathcal{V}_s})\bs{b}).
  \label{eq:def_S_k_correction}
\end{equation}
However, the spectral information contained in $\bs{x}_0$ can be removed from the correction equation~\eqref{eq:correction_equation}, suggesting the use of the Krylov space:
\begin{equation}
  \mathcal{S}_k=\mathcal{K}_k(\mat{P}_i\A,\mat{P}_j(\I-\mat{P}_{\mathcal{V}_s,\A\mathcal{V}_s})\bs{b}).
  \label{eq:def_S_k_projected}
\end{equation}
instead of the one of equation~\eqref{eq:def_S_k_correction}. In equation~\eqref{eq:def_S_k_projected}, the operators $\mat{P}_i,\mat{P}_j$ are projectors whose purpose is that of removing the information contained in the recycling space $\V_s$ from the correction equation~\eqref{eq:correction_equation}.

Some possible choices for $\mat{P}_i$ are:
\begin{itemize}
  \item $\mat{P}_{\mathcal{V}_s}:=\I-\V_s\V_s^T$: an orthogonal projector onto $\mathcal{V}_s$;
  \item $\mat{P}_{\mathcal{V}_s,\A\mathcal{V}_s}:=\I-\A\V_s\E_s^{-1}\V_s^T$: an oblique projector onto $\V_s$, orthogonal to $\A\mathcal{V}_s$;
  \item $\mat{P}_{\A\mathcal{V}_s,\A\mathcal{V}_s}:=\I-\A\V_s(\V_s^T\A^T\A\V_s)^{-1}\V_s^T\A^T$: an orthogonal projector onto $\A\V_s$;
\end{itemize}
while for $\mat{P}_j$ it is possible to choose $\mat{P}_j=\I$, which corresponds to augmentation methods, or $\mat{P}_j=\mat{P}_i$ for deflation methods.

An advantage of the projector $\mat{P}_{\mathcal{V}_s}$ is that its application to a vector requires fewer operations than $\mat{P}_{\mathcal{V}_s,A\mathcal{V}_s}$. However, $\mat{P}_{\mathcal{V}_s}$ has maximum rank $n-s$ only if $\V_s^T\V_s=\I_s$, while $\mat{P}_{\mathcal{V}_s,A\mathcal{V}_s}$ has maximum rank any time the columns of $\V_s$ are linearly independent. A more general variant of $\mat{P}_{\mathcal{V}_s}$ is:
\begin{equation}
  \mat{P}_{\mathcal{V}_s}:=\I-\V_s(\V_s^T\V_s)^{-1}\V_s^T,
  \label{eq:P_A_variant}
\end{equation}
that works even if $\V_s$ is not orthogonal.

Regarding the third choice, $\mat{P}_{\A\mathcal{V}_s,\A\mathcal{V}_s}$, we remark that the product $\A^T\A$ may increase the condition number of the problem, as is well-known for methods such as CGNE or CGLS~\cite{Saad}.

In the following, we will consider the six recycling methods of table~\ref{tab:which_methods}. All of these methods can be implemented by constructing a (left) preconditioner $\M$, and by considering a specific starting vector $x_0$. For example, the \emph{deflated oblique} method is obtained by considering $0$ as the starting vector, and $\M = (\I-\A\V_s\E_s^{-1}\V_s^T)$ as the (left) preconditioner, while the \emph{deflated orthogonal} method is obtained by using the (left) preconditioner $\M =(\I-\V_s\V_s^T)$ and $\bs{x}_0 = \A\V_s\E_s^{-1}\V_s^T \bs{b}$.

Lastly, we remark that the \emph{augmented oblique} and \emph{deflated oblique} methods are equivalent in the non-preconditioned case, since the initial vector for the latter method would be:
\begin{equation}
  (\I-\A\V_s\E_s^{-1}\V_s^T)^2\bs{b}=(\I-\A\V_s\E_s^{-1}\V_s^T)\bs{b}=\bs{r}_0.
  \label{eq:equivalence_oblique}
\end{equation}
However, since in general the presence of a preconditioner would prevent the first equality in equation~\eqref{eq:equivalence_oblique} from holding true, we consider the two cases as independent.

Regarding the methods of Table~\ref{tab:which_methods}, we remark that the \emph{deflated oblique} method coincides with the class of preconditioners studied in~\cite{NabbenVuik2006} if $\V$ contains an approximation of some eigenvectors of $\A$. The \emph{deflated LS} method is equivalent to the inner iteration of the GCRO method introduced in~\cite{deSturler96}. The other methods of Table~\ref{tab:which_methods} are, to the best of our knowledge, new.
Finally, we mention that alternative, but not equivalent classifications of augmented and deflated methods are available in the references~\cite{Gaul2013,Gutknecht2014}.

\begin{table}
  \caption{Krylov space used by the recycling methods considered in this work.}
\label{tab:which_methods}
\centering
\begin{tabular}{lrr}
\toprule
acceleration method        & Krylov operator & Krylov initial vector \\
\midrule
augmented orthogonal & $(\I-\V_s\V_s^T)\A$ & $\bs{r}_0$ \\
augmented oblique    & $(\I-\A\V_s\E_s^{-1}\V_s^T)\A$ & $\bs{r}_0$ \\
augmented LS    & $(\I-\A\V_s(\V_s^T\A^T\A\V_s)^{-1}\V_s^T\A^T)\A$ & $\bs{r}_0$ \\
deflated orthogonal  & $(\I-\V_s\V_s^T)\A$ & $(\I-\V_s\V_s^T)\bs{r}_0$ \\
deflated oblique     & $(\I-\A\V_s\E_s^{-1}\V_s^T)\A$ & $(\I-\A\V_s\E_s^{-1}\V_s^T)\bs{r}_0$ \\
deflated LS    & $(\I-\A\V_s(\V_s^T\A^T\A\V_s)^{-1}\V_s^T\A^T)\A$ & $(\I-\A\V_s(\V_s^T\A^T\A\V_s)^{-1}\V_s^T\A^T)\bs{r}_0$ \\
\bottomrule
\end{tabular}
\end{table}

\section{Acceleration by initial guess selection}
\label{sec:accel-init-guess}

The discussion so far has been quite general, but in the context of sequences of linear systems such those of equation~\eqref{eq:sequence_linear_systems}, arising e.g. from the discretization of a time-dependent PDE, it is reasonable to expect that some information could be extracted from the solutions of previous linear systems, and successfully used to accelerate the solution of successive linear systems.

For time-dependent problems, it is often standard to use an extrapolation based on a few preceding steps as the initial guess for the current linear system, on the basis that if the timestep is sufficiently small, the solutions will not differ too much from each other.

Following this basic idea, we go one step further and explore a strategy based on \emph{projection} rather than extrapolation. These methods can be interpreted as recycling methods, if the projection step is done coherently with the Krylov solver used in the solution of the linear system.

Recycling of Krylov subspaces for sequences of linear systems with changing matrices and right-hand sides was proposed in~\cite{deSturler2006}, where a generalization of the GCROT~\cite{deSturler99} and GMRES-DR~\cite{Morgan2002} methods, and a new GCRO-DR method are applied to Hermitian and non-Hermitian problems.
The effectiveness of the recycling methods mentioned above was proven in practical applications in~\cite{Kilmer2006,Wang2007,Mello2010}. Recycling BICG and BICGStab methods were developed in~\cite{Ahuja2012} and~\cite{Ahuja2015} respectively, where the recycling techniques play a key role in the set up of very efficient model order reduction methods. It is also possible to combine deflation with augmentation, as shown in~\cite{Gaul2013}, but this combined approach is not considered here.

The selection of a good initial guess for a Krylov space method has been studied especially for symmetric positive definite matrices not depending on time, with multiple right-hand sides~\cite{Fischer98,Erhel2000}.

For the model problem of Equation~\eqref{eq:sequence_linear_systems}, we are interested in providing good initial guesses for sequences of nonsymmetric, linear systems where both the matrix and the right-hand side change at each time step (but we suppose that the dimension of the linear systems remains costant).
Furthermore, in later sections we develop on this idea by selecting only a subset of previous solutions.

\subsection{Extrapolation}
\label{sec:extrapolation_methods}

Extrapolation methods are used to construct an initial guess for the iterative solver by Lagrangian interpolation from previous solutions. 
An initial guess $\bs{x}^i_0=\bs{x}_0(t_i)$ for the current linear system $\A_i\bs{x}^i=\bs{b}^i$ can be derived imposing that the $\bs{x}^i_0$ lie on the family of planes passing through $\bs{x}^{i-1}$ and $\bs{x}^{i-2}$. The expression for $\bs{x}^i_0$ can be obtained formally by introducing a pseudo time step $h$ and requiring that the function
\[
  \bs{f}(\tau)=c_1\bs{x}^{i-2}\tau + c_2\bs{x}^{i-1}\tau + c_3
\]
satisfies $\bs{f}(-2h)=\bs{x}^{i-2}$, $\bs{f}(-h)=\bs{x}^{i-1}$, and then:
\begin{equation}
  \bs{x}^i_0=\bs{f}(0)=2\bs{x}^{i-1}-\bs{x}^{i-2}.
  \label{eq:ext2}
\end{equation}
Similar expressions can be obtained using quadratic interpolation formulas:
\begin{equation}
  \bs{x}^i_0=3\bs{x}^{i-1}-3\bs{x}^{i-2}+\bs{x}^{i-3} 
  \label{eq:ext3}
\end{equation}
and cubic interpolation formulas:
\begin{equation}
  \bs{x}^i_0=4\bs{x}^{i-1}-6\bs{x}^{i-2}+4\bs{x}^{i-3}-\bs{x}^{i-4}. 
  \label{eq:ext4}
\end{equation}
Usually there is no advantage in further increasing the order of the interpolation polynomial, due to instabilities of equispaced interpolation for high order polynomials. Furthermore, the desired initial guess lies outside the interpolation interval, hence the name \emph{extrapolation}.

The full expression of the $k$-th iterate for the linear extrapolation method becomes:
\begin{equation}
  \bs{x}^i_k = 2\bs{x}^{i-1}-\bs{x}^{i-2} + \mat{S}_k(\C_k^T\A_i\mat{S}_k)^{-1}\C_k^T(\bs{b}^i-2\A_i\bs{x}^{i-1}+\A_i\bs{x}^{i-2}),
  \label{eq:k_iteration_ext}
\end{equation}
showing that the projector $\mat{P}_{\A\mathcal{S}_k,\mathcal{C}_k}$ is not modified by the extrapolation method, whose only effect is to perturb the right-hand side so that the special form of the initial guess is taken into account. In practice, nothing changes from the standard case of equation~\eqref{eq:x_k_expression}.

The main advantages of this method are its simplicity and low storage and work requirements. To the best of our knowledge, there exist no theoretical estimates on the effectiveness of extrapolation methods in reducing the initial residual.

\subsection{Simple projection}
\label{sec:projection_methods}

The idea of constructing a good initial guess through the extrapolation method of section~\ref{sec:extrapolation_methods} can be greatly improved if the interpolation on a set of previous solutions is replaced by projection.
This approach was first introduced in~\cite{Ali2012} for time-dependent linear systems with a fixed matrix.
The analysis we carry over differs from the original work~\cite{Ali2012}, leading us to introduce a new initial guess selection based on an instance of the Multiple Subspace Correction Method~\cite{XuZikatanov}.

Suppose that the $s$ previous solutions $\{\bs{x}^{i-j}\}_{j=1}^s$ are orthonormalized, e.g. through a modified Gram--Schmidt algorithm and stored columnwise in the matrix $\V_s$. Then we can look for an initial guess in the range of $\V_s$, that we name $\mathcal{V}_s=\Range{\V_s}$, $\bs{x}^k_0\in \V_s$, by requiring that the residual is orthogonal to the column space of $\V_s$:
\begin{equation}
  \bs{x}_0 = \V_s \bs{t}_s, \quad \V_s^T \bs{r}_0 = 0 \quad \Longrightarrow \quad \V_s^T\A\V_s\bs{t}_s=\V_s^T\bs{b},
  \label{eq:projection_guess}
\end{equation}
so that 
\begin{equation}
  \bs{x}_0 = \V_s\bs{t}_s=\V_s(\V_s^T\A\V_s)^{-1}\V_s^T\bs{b}=\V_s\E_s^{-1}\V_s^T\bs{b}.
  \label{eq:projection_t_s}
\end{equation}
A standard iterative method of choice can then be applied to the full linear system using the initial guess computed as above. In this case the iterations satisfy:
\begin{equation}
  \bs{x}_k=\V_s\E_s^{-1}\V_s^T\bs{b}+\mat{S}_k(\C_k^T\A\mat{S}_k)^{-1}\C_k^T(I-\A\V_s\E_s^{-1}\V_s^T)\bs{b}.
  \label{eq:iteration_projection}
\end{equation}

From equation~\eqref{eq:iteration_projection}, we see that the solution at the $k$-th iterate is obtained as a projection of the right-hand side on a suitable sequence of spaces.

Introducing the projector on $\mathcal{V}_s=\Range{\V_s}$, orthogonal to $\A\mathcal{V}_s$, i.e.,  $\mat{P}_{\mathcal{V}_s,\A\mathcal{V}_s}:=\A\V_s\E^{-1}\V_s^T$, iterate~\eqref{eq:iteration_projection} can be rewritten as:
\begin{equation}
  \bs{x}_k = \A^{-1}\left(\mat{P}_{\mathcal{V}_s,\A\mathcal{V}_s} + \mat{P}_{\A\mathcal{S}_k,\mathcal{C}_k}(\I-\mat{P}_{\mathcal{V}_s,\A\mathcal{V}_s})\right)\bs{b},
  \label{eq:two_projections}
\end{equation}
where it is clear that the correction of this method is based on a sequence of two projectors.
However, we cannot conclude that this method is an instance of the Multiple Subspace Correction Method because in general $\mathcal{V}_s$ is not orthogonal to the Krylov space $\A\mathcal{S}_k$, i.e., $\mathcal{V}_s\cap\A\mathcal{S}_k\ne\{\bs{0}\}$.

\subsection{Multiple Subspace Correction}
\label{sec:augmented_krylov}

To exploit the best approximation properties of the Multiple Subspace Correction Method, the spaces $\mathcal{V}_s$ and $\mathcal{S}_k$ should be orthogonal, and this can be achieved by replacing the matrix $\A$ in the linear systems~\eqref{eq:sequence_linear_systems} with $\tilde \A := (\I-\mat{P}_{\mathcal{V}_s})\A$, where $\mat{P}_{\mathcal{V}_s}=\V_s\V_s^T$ is the orthogonal projector onto $\V_s$.

In this case, the iteration satisfies:
\begin{equation}
  \bs{x}_k^i=\V_s\E_s^{-1}\V_s^T\bs{b}+\mat{S}_k(\C_k^T(I-\V_s\V_s^T)\A\mat{S}_k)^{-1}\C_k^T(\I-\A\V_s\E_s^{-1}\V_s^T)\bs{b},
  \label{eq:augmented_iterate}
\end{equation}
a slight modification of the operator in equation~\eqref{eq:iteration_projection}. Also, the projector expression~\eqref{eq:two_projections} holds with the difference that in this case $\tilde\A\mathcal{S}_k$ and $\mathcal{V}_s$ are orthogonal. The higher computational costs due to the extra matrix-vector products required by the application of $\I-\mat{P}_{\mathcal{V}_s}$ are hopefully recovered by the best approximation estimate on the residual~(following the lines of~\cite{EiermannErnstSchneider}):
\begin{equation}
  \|\bs{y}-(\mat{P}_{\mathcal{V}_s,\A\mathcal{V}_s}\bs{y}+\mat{P}_{\A\mathcal{S}_k,\mathcal{C}_k}(\I-\mat{P}_{\mathcal{V}_s})\bs{y})\| = \min_{\bs{t}\in\mathcal{V}_s+\mathcal{S}_k}\|\bs{y}-\bs{t}\|.
  \label{eq:best_approximation_augmented}
\end{equation}
This property can be proved by noting that
\begin{equation}
  \mat{P}_{\mathcal{V}_s,\A\mathcal{V}_s}\bs{x}+\mat{P}_{\A\mathcal{S}_k,\mathcal{C}_k}(\I-\mat{P}_{\mathcal{V}_s})\bs{x}=\mat{P}_{\mathcal{V}_s,\A\mathcal{V}_s}\bs{x}+\mat{P}_{\A\mathcal{S}_k,\mathcal{C}_k}\bs{x}=\mat{P}_{\mathcal{V}_s+\mathcal{S}_k}\bs{x} \qquad \forall x\in\mathbb{R}^n,
  \label{eq:proof_min_projection}
\end{equation}
where $\mat{P}_{\mathcal{S}_k}\mat{P}_{\mathcal{V}_s}=0$ was used since $\mathcal{V}_s$ and $\mathcal{S}_k$ are orthogonal.

However, in many cases the spaces $\mathcal{V}_s$ and $\mathcal{S}_k$ can become relatively large, and the advantages of having larger approximation spaces may be invalidated by excessive computational and storage costs. To mitigate this effect, restarting and truncation methods have been introduced, meaning that the approximation spaces are not allowed to exceed a fixed dimension.

To conclude, we summarize the main advantages of using multiple projections over extrapolation methods:
\begin{itemize}
  \item the stability issues of equispaced interpolation are removed; hence it can be worthwhile to keep a larger number of previous solutions;
  \item for the projection method (Galerkin) best approximation-type estimates hold;
  \item the vector space for the projection step can be obtained in any way, for instance the previous solutions can be orthonormalized and augmented with any other kind of basis functions, regardless of their origin.
\end{itemize}

This last point is worth noting, since in general we are not interested in the previous solutions \emph{per se}, but only in how representative they are as solutions to the current step. In the next section, we claim that a good way to keep as much information as possible from previous solutions and at the same time maintaining the size of the projection space small is to use the Singular Value Decomposition (SVD) of a fixed number of preceding solutions, and only retain those modes that contain information related to the singular values of largest or smallest magnitude.  This strategy is also very effective when combined with deflation preconditioners, as shown in the next section.


\section{Construction of the recycling spaces}
\label{sec:construction_augmentation}

In the previous sections, we implied more or less explicitly that the augmentation spaces should be built by performing some orthonormalization procedure on a fixed number of previous solutions. In fact, there is complete freedom in the choice of the augmentation spaces, and the only obstructions are considerations on the approximation, on the algebraic stability, and on the effectiveness of the resulting algorithm.
In Section~\ref{sec:ritz-harmonic-ritz}, we shortly review some common methods for the selection of the augmentation spaces and the restarting techniques used when their dimension becomes too large, while in Section~\ref{sec:simplified_svd_truncation} we propose an alternative strategy, based on the construction of the enriching space using Singular Value Decomposition. 
This idea has been traditionally used for the construction of efficient Model Order Reduction methods~\cite{Kressner2011}, but its application to the construction of recycling spaces could be convenient over the computation of Ritz vectors, at least in the case with time-dependent matrices.

\subsection{Ritz and harmonic Ritz pairs}
\label{sec:ritz-harmonic-ritz}

The discussion of section~\ref{sec:adapted_deflation} suggests constructing the recycling spaces by computing the $s$ eigenvalues of $\A$ with the smallest magnitude and storing the respective eigenvectors column wise in the matrix $\V_s$. However, when $\A$ is too large to resort to algorithms such as QR subspace iterations for the approximation of the eigenvalues of $\A$, a few interesting eigenpairs can be approximated using the information generated by a GMRES cycle.

This approach was introduced by~\cite{Morgan}, see also~\cite{Morgan2002}, where it is shown that keeping approximate eigenvectors corresponding to small eigenvalues can greatly improve the convergence rate of the GMRES method after restart. Deflative preconditioners based on approximate eigenvalues of the linear operators are discussed also in~\cite{Kharchenko95} for the restarted GMRES method.
One further step in this direction was done in~\cite{Giraud2010}, with the introduction of a flexible GMRES with a deflation space built on some approximate eigenvectors, so that a different preconditioner could be used for each GMRES iteration.

Specifically, if $\V_k\in\mathbb{R}^{n\times k}$ is a rectangular matrix spanning the subspace $\mathcal{V}_s\subset\mathbb{R}^n$, the couple $(\theta,\bs{y})\in\mathbb{C}\times \V_s$ is called a \emph{Ritz pair}~\cite{PaigeParlettvdVorst,SaadEigen} of $\A$ with respect to the space $\mathcal{V}_s$ if:
\begin{equation}
  \bs{v}^H(\A\bs{y}-\theta\bs{y}) = 0\qquad\forall \bs{v}\in \mathcal{V}_s.
  \label{eq:def_Ritz_pair}
\end{equation}
Equivalently, $(\theta,\bs{y})$ is a Ritz pair for $\A$ with respect to $\mathcal{V}_s$ if $(\A-\theta \I)\bs{y}\in\mathcal{V}_s^\perp$. In this case, $\theta$ and $\bs{y}$ are called respectively \emph{Ritz value} and \emph{Ritz vector} for $\A$.
The Ritz pairs can be computed directly by restricting $\A$ to $\mathcal{V}_s$, and then computing the eigenvalues of the resulting Hessenberg matrix $\mat{H}_s=\V_s^T\A\V_s$:
\begin{equation}
  \mat{H}_s\bs{z}_i-\theta_i\V_s^H\V_s\bs{z}_i=\bs{0}\qquad \bs{y}_i=\V_s\bs{z}_i.
  \label{eq:eigenvalues_hessenberg}
\end{equation}
As the dimension of $\mathcal{V}_s$ is increased, the Ritz values converge (see e.g.~\cite[Ch.~34]{TrefethenNumLinAlg} and~\cite[Ch.~7]{Demmel}) to the external eigenvalues of $A$.

Alternatively, Ritz pairs can be replaced or enhanced through the computation of \emph{harmonic Ritz pairs}~\cite{PaigeParlettvdVorst,Sorensen}, defined as the solutions to the Petrov--Galerkin projection:
\begin{equation}
  \V_s^H\A^H\A\V_s\bs{z}_i-\theta_i\V_s^H\A^H\V_s\bs{z}_i=\bs{0}\qquad \bs{y}_i=\V_s\bs{z}_i.
  \label{eq:eigenvalues_harmonic}
\end{equation}

The harmonic Ritz values converge to the inner eigenvalues of $A$ with minimum distance from the origin of the complex plane.
Indeed, introducing the matrix $W_s=A\V_s$, problem~\eqref{eq:eigenvalues_harmonic} can be rewritten as:
\begin{equation}
  \W_s^H\W_s\bs{z}_i-\theta_i \W_s^H\A^{-1}\W_s\bs{z}_i=\bs{0}\qquad \bs{y}_i=\W_s\bs{z}_i,
  \label{eq:harmonic_ritz_inverse}
\end{equation}
showing that the harmonic Ritz values of $\A$ with respect to $\mathcal{V}_s$ coincide with the Ritz values of $\A^{-1}$ with respect to $\mathcal{W}_s=\A\mathcal{V}_s$.

\subsection{Simplified SVD truncation}
\label{sec:simplified_svd_truncation}
In this work we explore an alternative enriching strategy, that aims to retain the relevant information contained in the previous solution steps.  In particular, we collect a number of previous solution steps, perform a Singular Value Decomposition, and keep a fixed number of singular vectors associated with the higher singular values, that are then used as basis for the enriching subspace.
A similar strategy was used in a different context, namely model order reduction of parametrized matrix-valued linear systems, in~\cite{Kressner2011}.

The $m$ previous solutions $\bs{x}_m$ are stored in the matrix $\X_m\in\mathbb{R}^{n\times m}$. By computing an SVD of $\X_m$ we keep the $s$ left singular vectors with the largest (or smallest) singular values.


If we let $\X_m=\Psi_m\Sigma_m\Phi_m^T$ be the SVD of $\X_m$, and $\Psi_s$ the matrix whose columns are the $s$ largest left singular vectors, then the Schmidt--Eckart--Young theorem guarantees the following optimality property of the SVD: if $Y$ is a given matrix, the best approximation of $Y$ by a matrix of rank $s$ is given by:
  \begin{equation}
    Y_s := \sum_{i=1}^s \sigma_i \bs{\psi}_i\bs{\phi}_i^T,
    \label{eq:best_rank_s}
  \end{equation}
  where $Y=\Psi\Sigma\Phi^T$ is the SVD of $Y$, $\sigma_i$ is the entry in the $i$-th row of $\Sigma$, $\bs{\psi}_i$ and $\bs{\phi}_i$ the columns of $\Psi$ and $\Phi$ respectively.
  The optimality property given by the Schmidt--Eckart--Young theorem can be written explicitly as:
\begin{equation}
  Y_s=\arg\min_{\overset{\Xi\in\mathbb{R}^{n\times s}}{\rank\Xi=s}}\|\Y-\Xi\|_F.
  \label{eq:opt_svd}
\end{equation}

To clarify how this alternative truncation strategy could be implemented in a solver for a sequence of related linear systems, we refer to Algorithm~\ref{alg:simplified_truncation}.

\begin{algorithm}[tbp]
  \SetKwInOut{Input}{Input}
  \SetKwInOut{Output}{Output}
  \Input{Three nonnegative integers $m,s,\ell$; $m$ is the number of retained previous solutions, $s\le m$ is the dimension of the recycling space, $\ell\ge1$ is the interval measured in timesteps between two consecutive SVDs.}
  \BlankLine
    \For{$i=1,\dots$}{
    solve $\A_i\bs{x}_i=\bs{b}_i$ with any recycling scheme of section~\ref{sec:synoption_aug_def}\;
    store the solution $\bs{x}_i$ in the matrix $\X_m\in\mathbb{R}^{n\times m}$: \[\X_m=[\bs{x}_{i-m},\dots,\bs{x}_i].\] If necessary, remove the $i-m-1$-th row from $\X_m$ \;
    \If{$(i\mod\ell)$=0}{
    compute the SVD decomposition of $\X_m$: $\X_m=\Psi_m\Sigma_m\Phi_m^T$\;
    update the recycling space $\V_s\in\mathbb{R}^{n\times s}$ by copying the first $s$ largest left singular modes of $\X_m$: \[\V_s=\Psi_m[1:s];\]
  }
  }
  \caption{The SVD-based truncation algorithm described in section~\ref{sec:simplified_svd_truncation}.}
  \label{alg:simplified_truncation}
\end{algorithm}

Note that some of the recycling algorithms presented in section~\ref{sec:synoption_aug_def} require the solution of a linear system involving the matrix $\E_s=\V_s^T\A\V_s$ as in equation~\eqref{eq:projection_guess}. The reader should be aware that in principle the matrix $\E_s$ may be singular when $\mathcal{S}_k$ and $\mathcal{C}_k$ are the same space.

However, for the choices of $\V_s$ used in this work, we never had to face such difficulties, hence we do not take any preventive measure against this breakdown possibility.

\section{Numerical experiments}
\label{sec:numerical_experiments}

In this section we compare the accelerating methods proposed in Sections~\ref{sec:synoption_aug_def} and~\ref{sec:accel-init-guess} on a standard test case for nonsymmetric matrices: a parabolic scalar convection-diffusion equation discretized with Finite Elements. 

\subsection{Scalar convection-diffusion}
\label{sec:convection_diffusion}

This benchmark consists of a scalar parabolic convection-diffusion equation on a time interval $[0,T]$ and a domain $\Omega\subset\mathbb{R}^2$, that in weak form reads~\cite{QuarteroniValli}: find $u\in L^2([0,T],H^1_0(\Omega))\cap C^0((0,T),L^2(\Omega))$ such that:
\begin{equation}
\begin{cases}
(v,\partial_t u)+\nu(\nabla v,\nabla u) - (v,u\bs{b}\cdot\nabla u)=(v,f) \qquad &\forall v\in H^1_0(\Omega) \\
u(t=0) = u_0, & u_0\in L^2(\Omega),
\end{cases}
\label{eq:convection_diffusion}
\end{equation}
where $\nu$ is the diffusion coefficient, $\bs{b}:\Omega\to\mathbb{R}^2$ is a given solenoidal vector field, and $f:(0,T)\times\Omega\to\mathbb{R}$ is a possibly time-dependent given forcing term. Here we choose $\Omega$ to be the unit square with coordinates $(x,y)\in[0,1]^2$ and $\bs{b}=(-\sin(\pi x)\cos(\pi y),\cos(\pi x)\sin(\pi y))$.

Following~\cite{BurgersTurb}, we model $f$ as acting randomly on the low frequency modes of the Laplace operator, namely:
\begin{equation}
f = \frac{C}{2}\sum_{j=1}^{16}c_j\exp\left(-\frac{j^2}{20}\right)\sin(2j\pi x)\sin(2j\pi y),
\end{equation}
where $C=0.1$, $c_1=1$ and $\{c_j\}_{j=2}^{16}\subset[-1,1]$ are random coefficients sampled at every timestep.
The rationale of this choice is to avoid the solution settling to a stationary configuration, so that the turbulent behaviour of one-dimensional Burgers equation can be simulated in this more complex two-dimensional setting.

Treating equation~\eqref{eq:convection_diffusion} with a finite difference method in time (in this case an implicit Euler method) and with a first order finite element method, we get the sequence of discrete problems: for $n\in[1,N_T]$, find $u_h^n\in V_h$ such that
\begin{equation}
\begin{cases}
  \frac{1}{\Delta_t}(v,u_h^n-u_h^{n-1})+\nu(\nabla v,\nabla u_h^n) + (v,u_h^{n-1}\mat{P}_{V_h}\bs{b}\cdot\nabla u_h^n)=(v,\mat{P}_{V_h}f) \qquad &\forall v\in V_h \\
u_h^0=\mat{P}_{V_h} u_0,
\end{cases}
\end{equation}
where $\Delta t$ is the time step for the finite difference discretization, and $N_T=T/\Delta_t$ is the number of timesteps.
For the finite element space $V_h$, we adopt first order Lagrange elements on a grid formed by $64\times64$ uniform square elements, as provided by the \texttt{dealii} library~\cite{dealII84}.

We design three test cases, with $\nu=10^{-1}$, $\nu=10^{-2}$ and $\nu=10^{-3}$. Values of $\nu$ lower than $10^{-3}$ are not stable in the grid described above.
For each test case, we run a total of 1000 timesteps of width $\Delta t=0.5$, with a zero initial value and with a Symmetric Successive Over-Relaxation (SSOR) preconditioner. The iterative solver is halted when the 2-norm of the residual relative to the 2-norm of the right-hand side is below a tolerance set to $10^{-8}$.

In the first test, we compare the extrapolation initial guess described in Section~\ref{sec:extrapolation_methods} with the projection initial guess of Section~\ref{sec:projection_methods}. This test is performed with $\nu=10^{-1}$, and the results are reported in Table~\ref{tab:iterations_ext_proj}, where it is possible to see how the projection-based initial guess reduces the number of GMRES iterations to converge. Furthermore, increasing the dimension of the recycling space does not increase the iteration count for the projection method, while it clearly worsens the convergence rate of the extrapolation method.
\begin{table}
  \caption{Average number of iterations (av) and standard deviation of the iteration count (dev) for the extrapolation (ext) and projection (proj) av initial guess, as a function of the recycling space dimension (s). The results are for the case $\nu=10^{-1}$.}
\label{tab:iterations_ext_proj}
\centering
\begin{tabular}{ccccc}
\toprule
s      & ext av & ext stddev & proj av & proj stddev \\
\midrule
2      & $65.31$ & $3.86$ & $62.38$ & $3.56$ \\
3      & $67.57$ & $4.23$ & $61.17$ & $3.39$ \\
4      & $70.17$ & $4.20$ & $60.72$ & $3.06$ \\
\bottomrule
\end{tabular}
\end{table}

After verifying the expected behaviour for the two prediction methods, we compare the 4 recycling solvers of Table~\ref{tab:which_methods}.

To get some insight on the influence of the recycling parameters on the convergence properties of the iterative methods, we set up a campaign of simulations as follows.
We consider the same problem described above for the case $\nu=10^{-2}$, and run the various recycling algorithms letting the dimension of the Krylov space\footnote{equal to the GMRES restart value.} and the number of saved solutions sweep between 10 and 60 (in steps of 2). The dimension of the recycling space is chosen equal to the GMRES restart. Each campaign of simulations requires $26^2=676$ runs.

\begin{figure}[tbp]
  \centering
  \includegraphics[width=0.50\columnwidth]{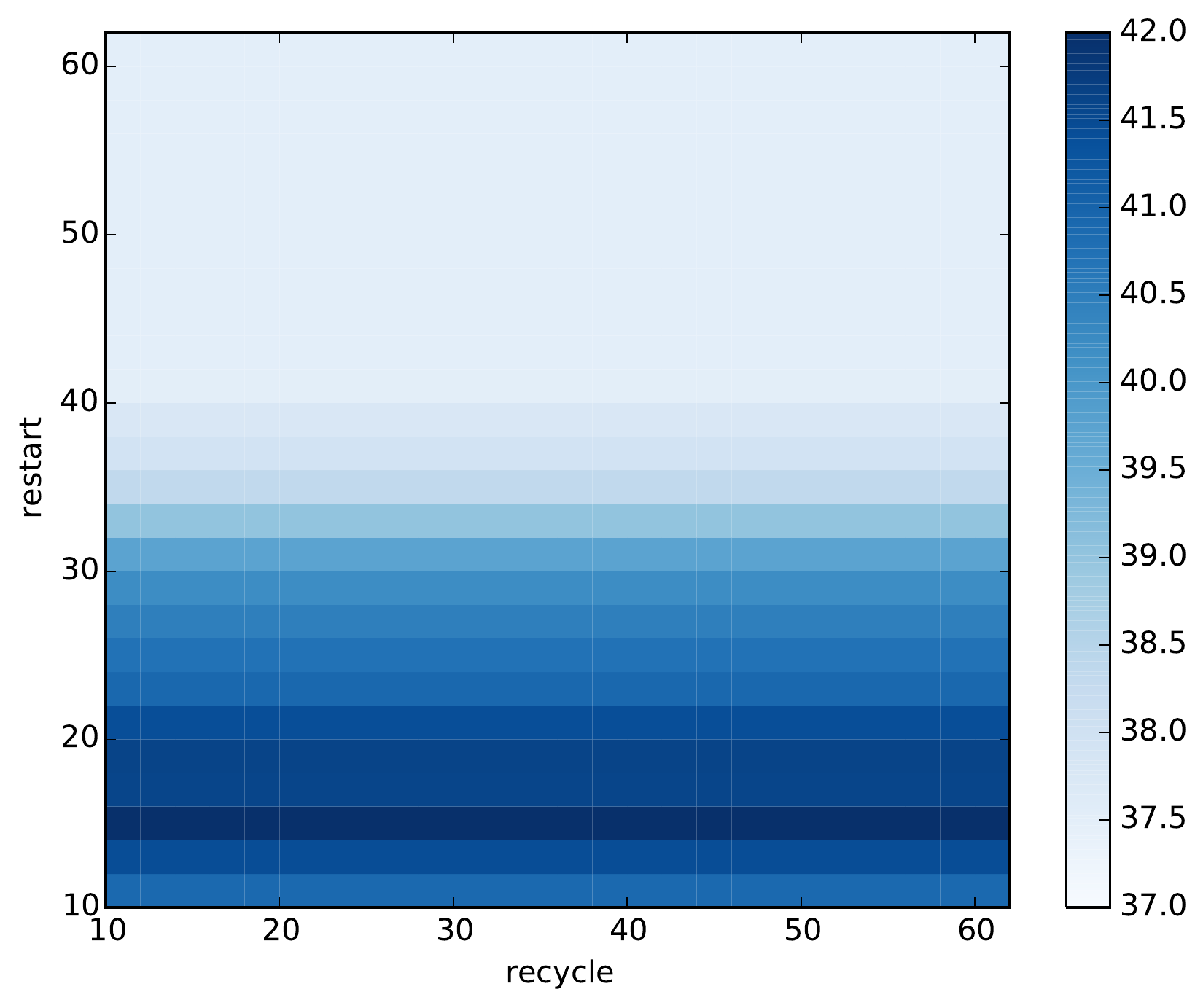}
  \caption{Average number of iterations to convergence without applying the recycling methods. Note that in this case the horizontal axis does not denote a true variable, since the recycling is not present. However, we keep the two-dimensional map style for consistency with the next figures.}
\label{fig:map_unp}
\end{figure}

At first, we run a standard solver without recycling, using only a symmetric successive over relaxation (SSOR) preconditioner. From the results, shown in Figure~\ref{fig:map_unp}, it can be seen that increasing the dimension of the Krylov space (higher values of the restart parameter) yields convergence in fewer iterations, although there are counterexamples to this behaviour as shown in~\cite{Embree}. In this case there is no variation of the iteration number in the horizontal direction, since the recycling method is not set up.

\begin{figure}[tbp]
  \centering
  \includegraphics[width=0.49\columnwidth]{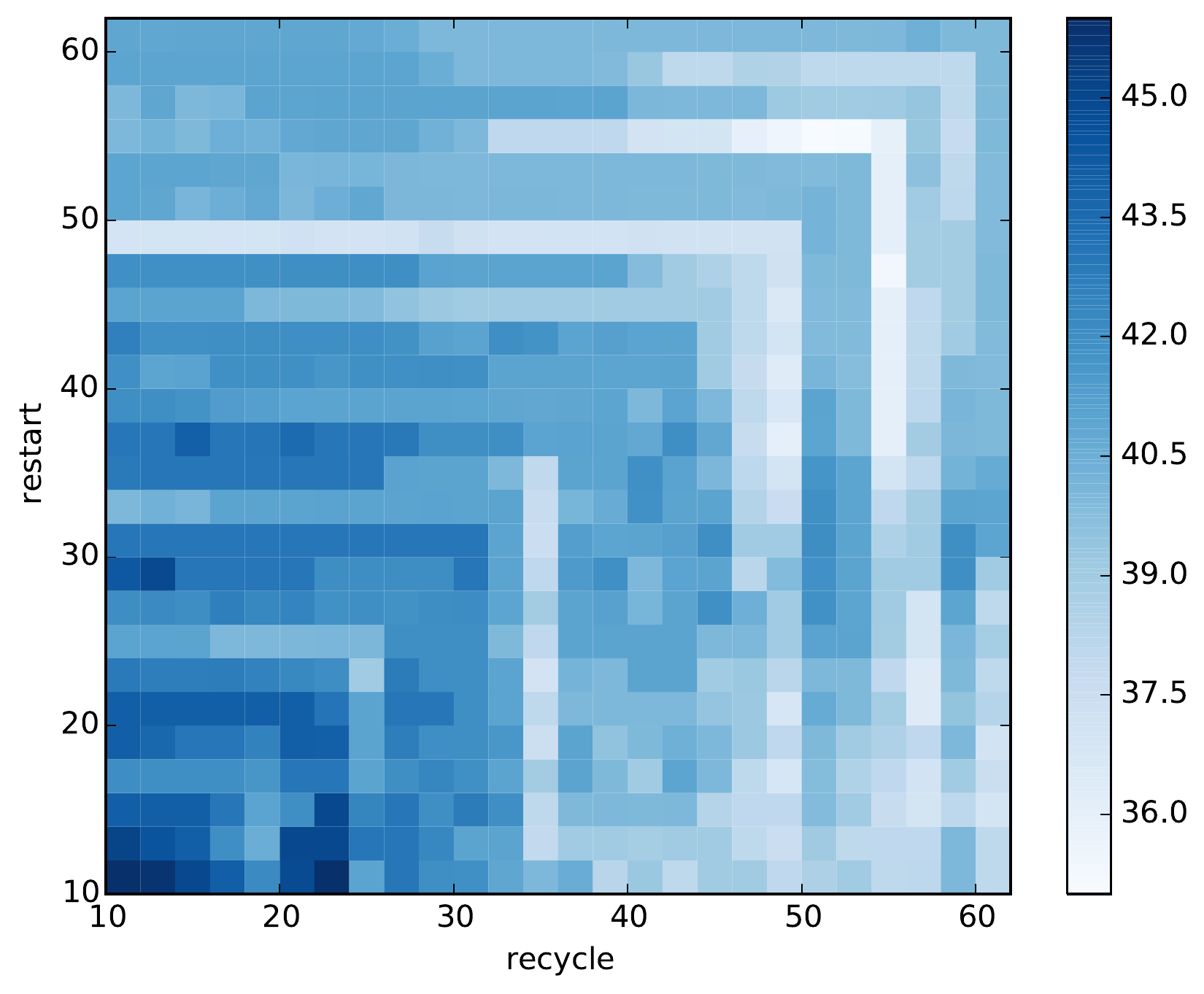}
  \includegraphics[width=0.49\columnwidth]{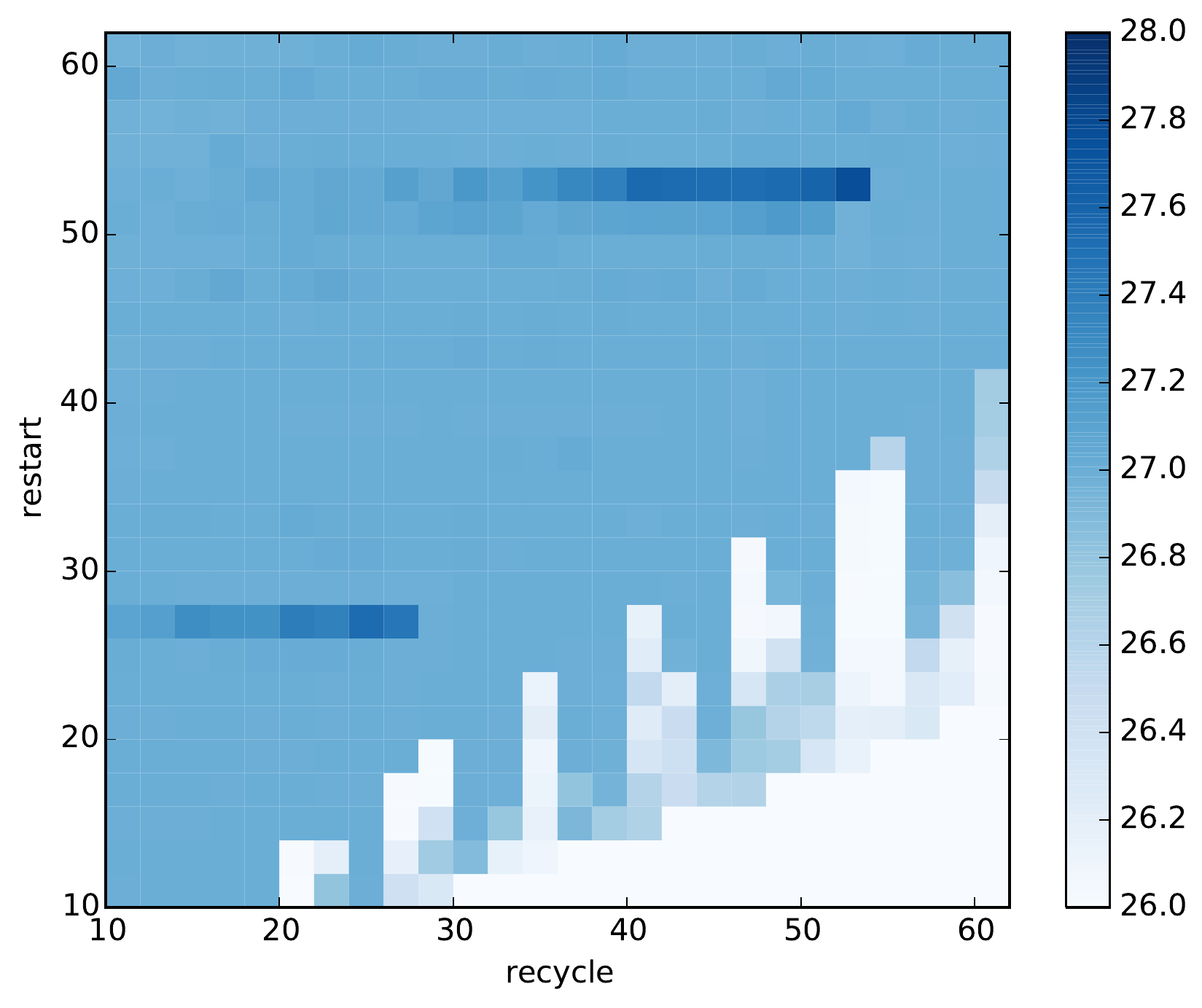}
  \caption{Average number of iterations to convergence for the \emph{augmented orthogonal} (left) and for the \emph{augmented oblique} (right) methods of table~\ref{tab:which_methods}. In both figures, the SVD selection is performed at the end of each time step, and the SSOR preconditioner is used.}
  \label{fig:map_aug}
\end{figure}

Then, we run the same campaign of simulations for the \emph{augmented orthogonal} and \emph{augmented oblique} recycling methods, corresponding to the first two rows of Table~\ref{tab:which_methods}. As in the previous case, the SSOR preconditioner is used within the recycling method, and the SVD selection is performed at every time step. As reported in Figure~\ref{fig:map_aug}, in most of the runs the augmentation methods are reducing the iteration count with respect to the unrecycled case of Figure~\ref{fig:map_unp}.
In both cases it is difficult to figure out a clear pattern, however in average we observe that increasing the dimension of the recycling space reduces the iteration count.

\begin{figure}[tbp]
  \centering
  \includegraphics[width=0.49\columnwidth]{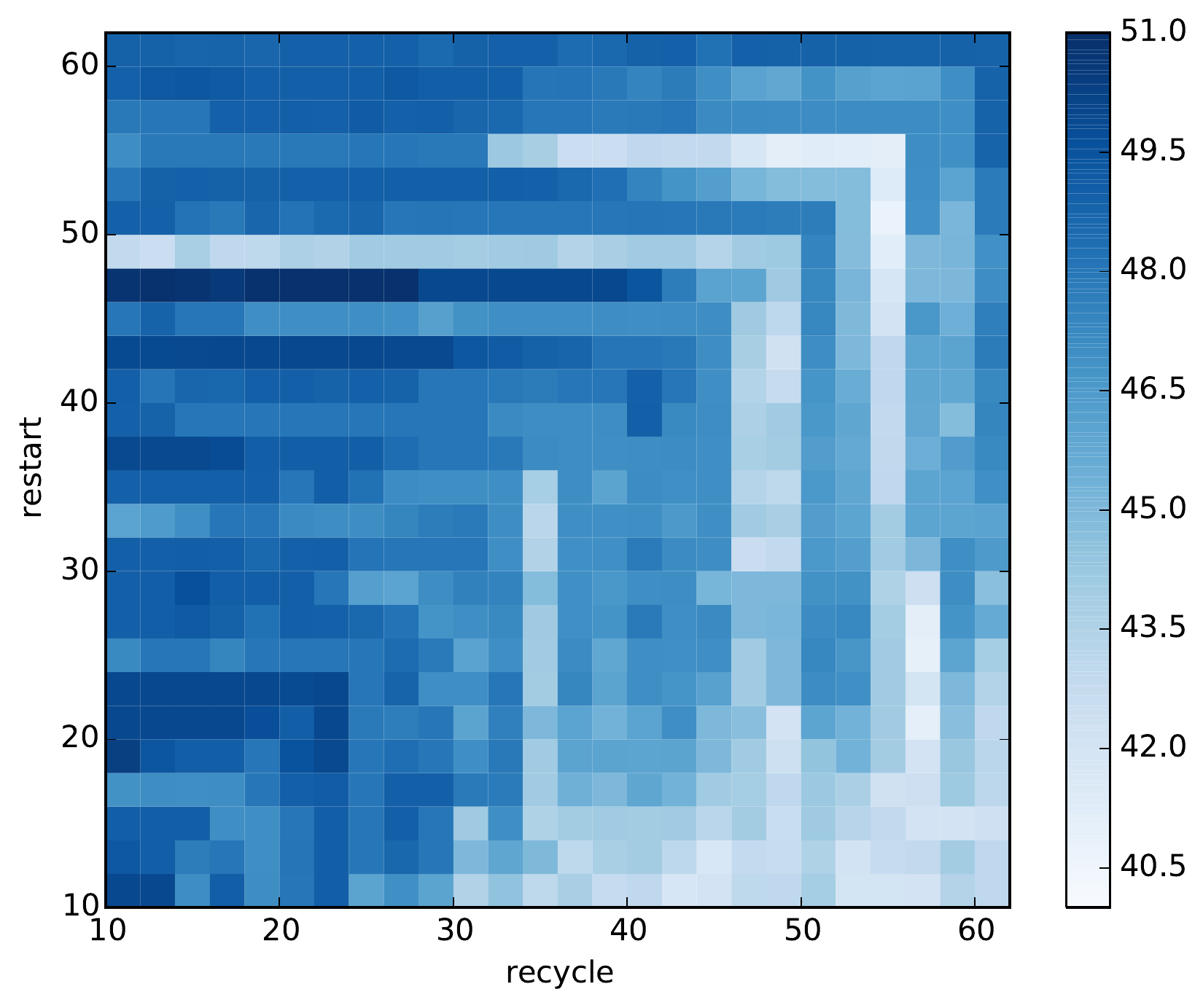}
  \includegraphics[width=0.49\columnwidth]{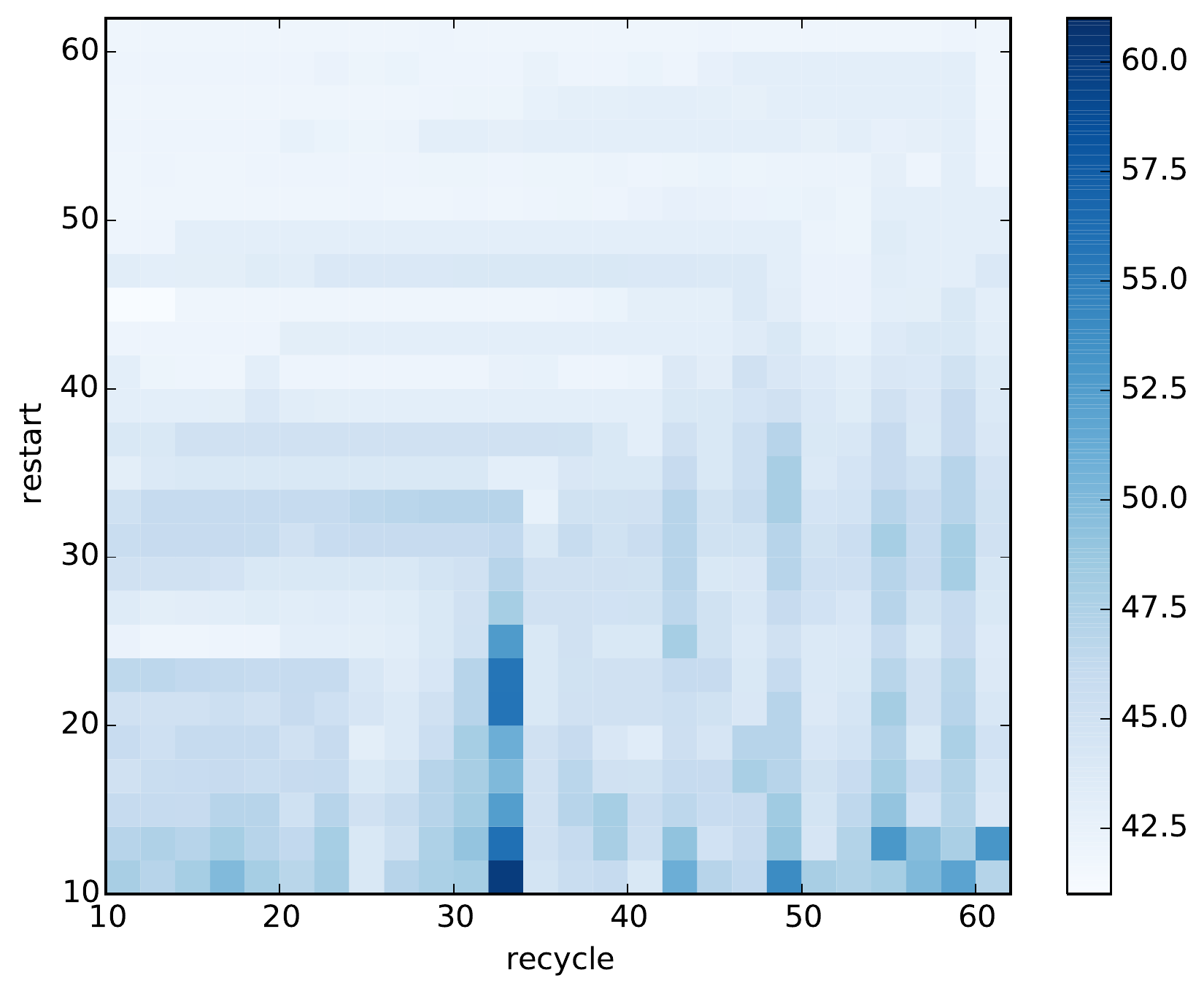}
  \caption{Average number of iterations to convergence for the \emph{deflated orthogonal} (left) and for the \emph{deflated oblique} (right) methods of table~\ref{tab:which_methods}. In both figures, the SVD selection is performed at each time step, and the SSOR preconditioner is used.}
  \label{fig:map_def}
\end{figure}

The second set of runs involves the \emph{deflated orthogonal} and \emph{deflated oblique} methods, corresponding to the fourth and fifth rows of Table~\ref{tab:which_methods}. As for the previous two cases, the numerical experiments take advantage of the SSOR preconditioner, and the SVD selection is performed at the end of each time step.
The results of this test are reported in Figure~\ref{fig:map_def}.
It can be seen that the \emph{augmented orthogonal} and the \emph{deflated orthogonal} methods lead to visually similar convergence patterns, corresponding to Figures~\ref{fig:map_aug} (left) and~\ref{fig:map_def} (left). However, the augmentation method appears to be superior, sparing a few iterations for almost all the cases.
By comparing Figures~\ref{fig:map_aug} (right) and~\ref{fig:map_def} (right), we see that the oblique augmentation method performs much better than its deflated counterpart, since in this case the improvement amounts to almost halving the iteration count.
In the map on the right of Figure~\ref{fig:map_def}, it is possible to recognise an interesting pattern relating faster convergence to higher values of the restart parameter.

\begin{figure}[tbp]
  \centering
  \includegraphics[width=0.49\columnwidth]{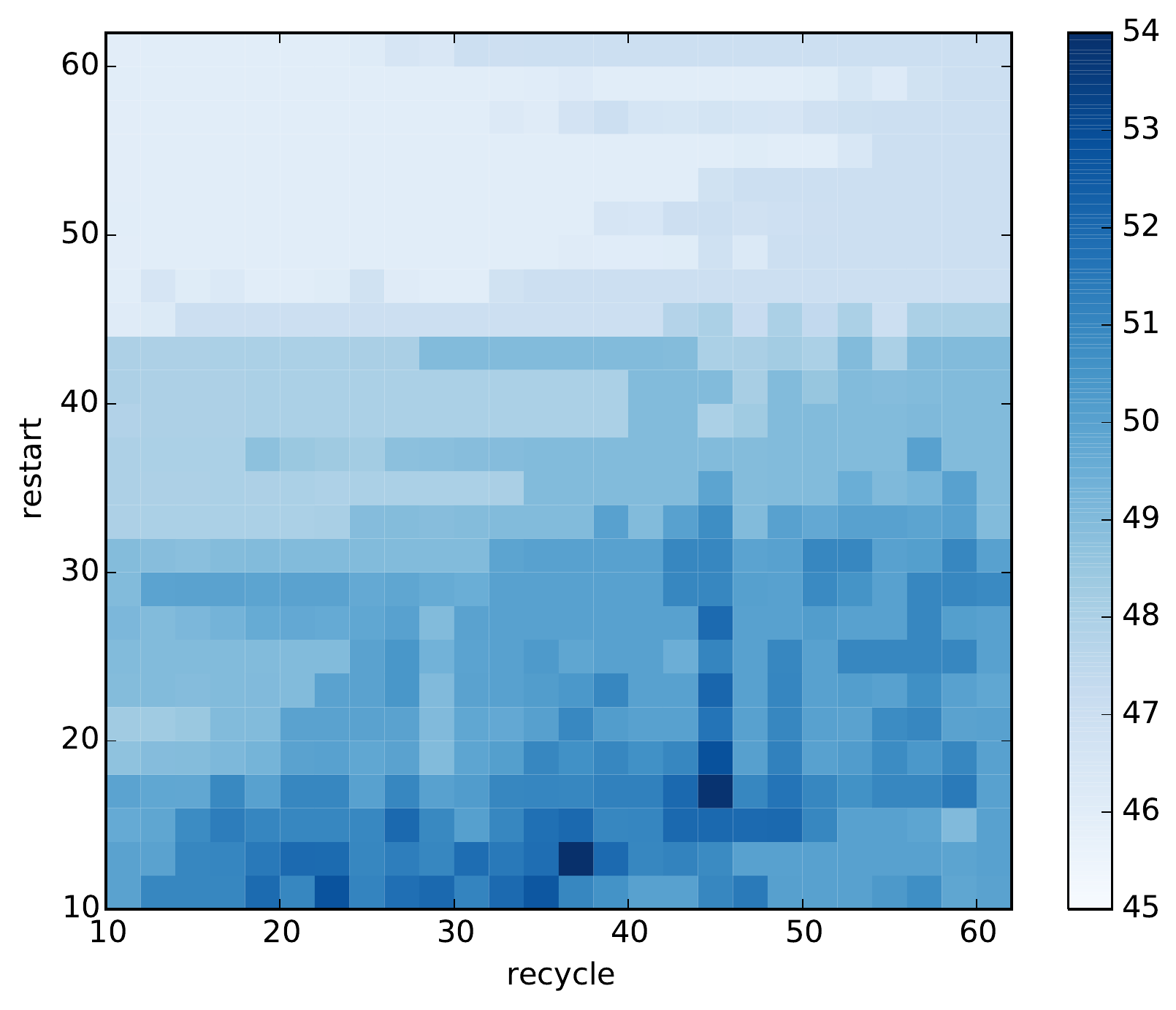}
  \caption{Average number of iterations to convergence for the \emph{deflated LS} (right) method of table~\ref{tab:which_methods}. The SVD selection is performed at each time step, and the SSOR preconditioner is used.}
  \label{fig:map_w}
\end{figure}

Next, we compare the performance of the two \emph{orthogonal LS} acceleration methods (rows 3 and 6 of Table~\ref{tab:which_methods}). We have not been able to make the \emph{agumented LS} method converge, so in Figure~\ref{fig:map_w} we show only the results for the \emph{deflated LS} method.
Overall, it seems that despite the higher computational cost due to the higher number of matrix-vector multiplications, this strategy is not paying back with an increased convergence rate.

\begin{figure}[tbp]
  \centering
  \includegraphics[width=0.49\columnwidth]{map_def_a.pdf}
  \includegraphics[width=0.49\columnwidth]{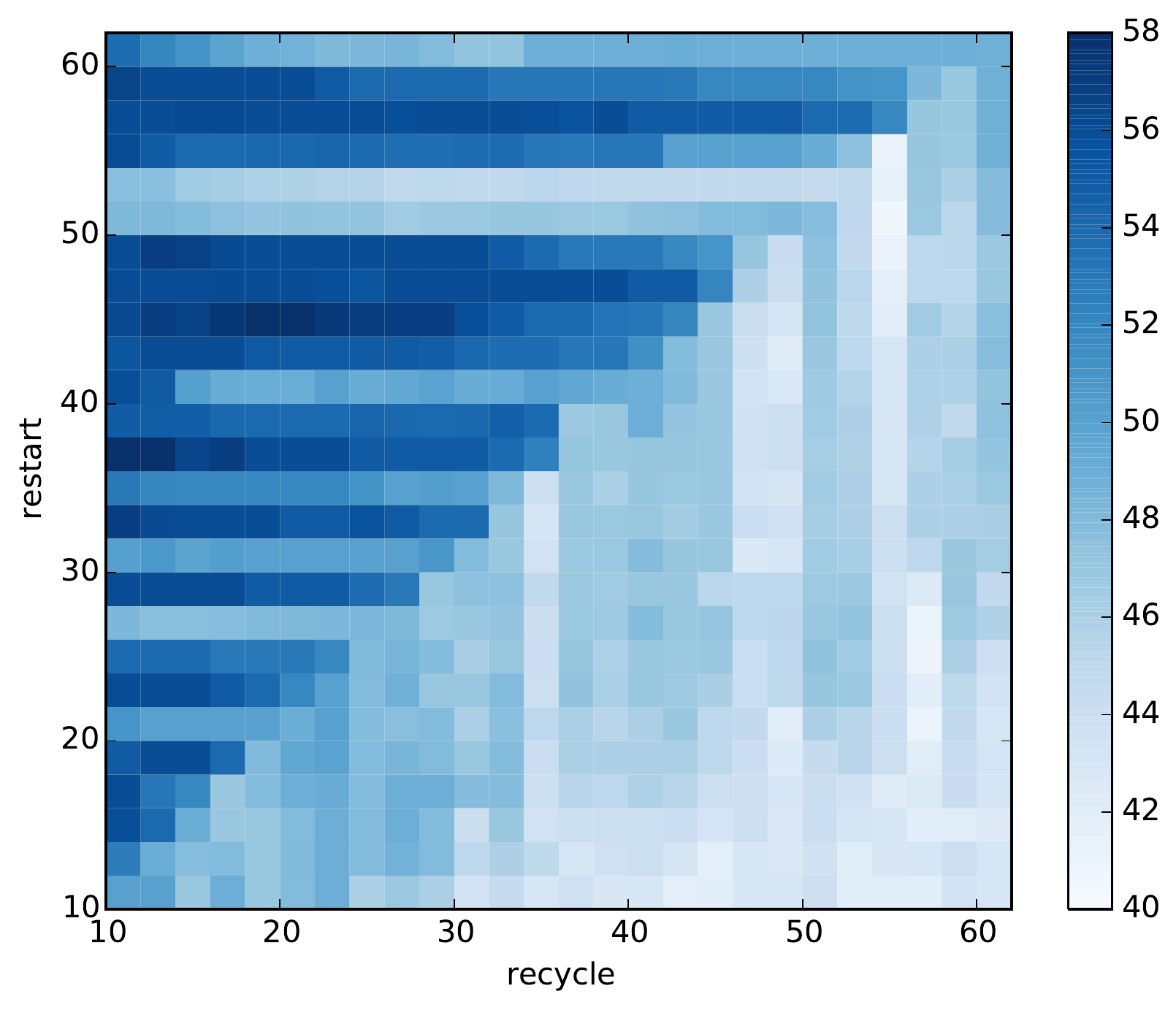}
  \includegraphics[width=0.49\columnwidth]{map_def_d.pdf}
  \includegraphics[width=0.49\columnwidth]{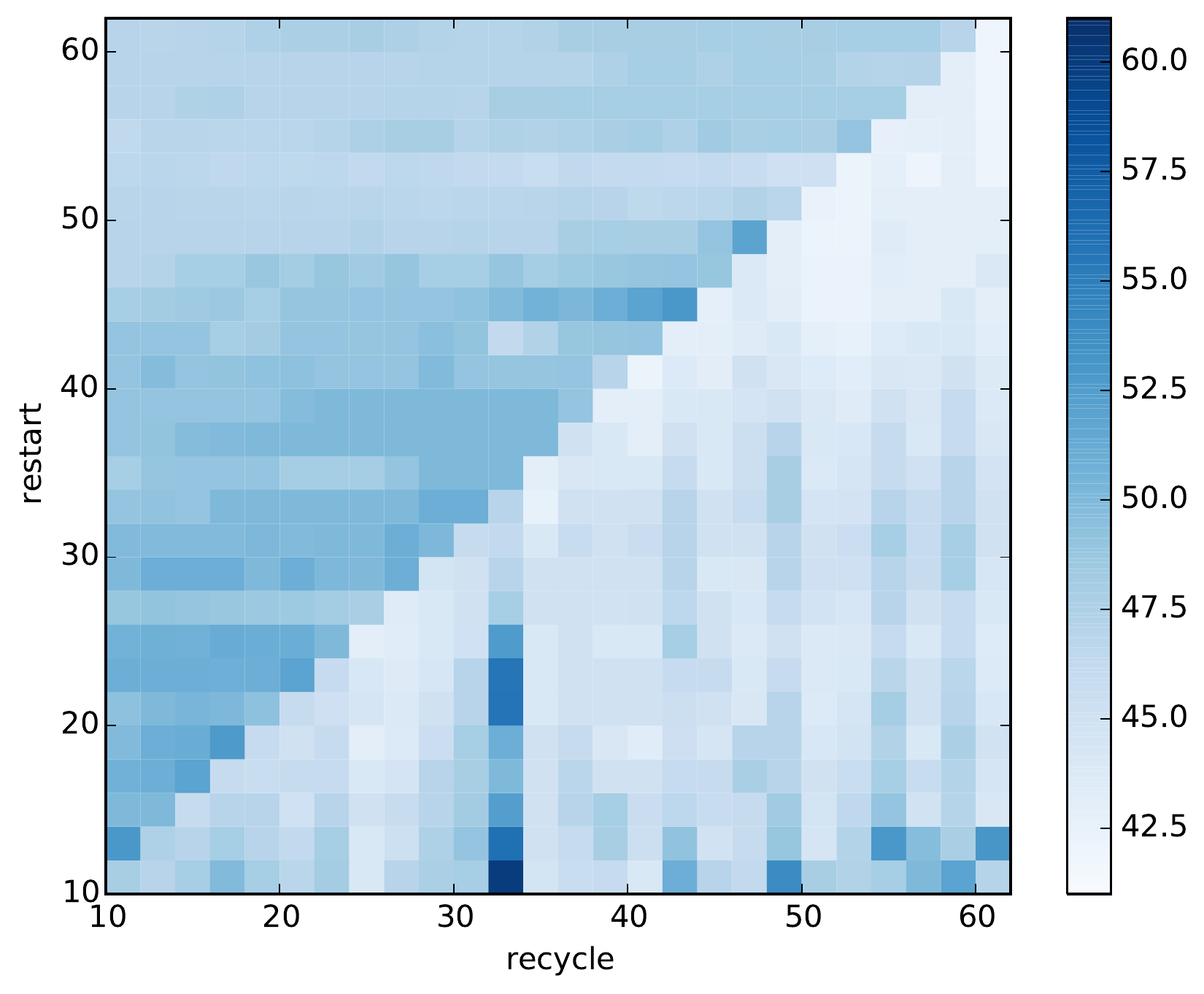} \\
  \caption{Average number of iterations to convergence for the \emph{deflated orthogonal} (top row) and for the \emph{deflated oblique} (bottom row) methods of table~\ref{tab:which_methods}. On the left column, the recycling subspace is spanned by the SVD modes with higher singular values, while on the right column the recycling subspace is formed by the POD modes with smaller singular values. In every case, the SVD selection is performed at each time step, and the SSOR preconditioner is used.}
  \label{fig:cmp_def_recycling_space}
\end{figure}

As a last test, we compare the two deflative preconditioners with different choices of the recycling subspace. The first choice consists in the SVD modes associated with higher singular values, the second choice conversely consists in the SVD modes with smaller singular values.
The comparison for the two deflative preconditioners is available in Figure~\ref{fig:cmp_def_recycling_space}.
For both the deflative preconditioners, we see that the second choice for the recycling subspace introduces a clear diagonal pattern in the diagrams of Figure~\ref{fig:cmp_def_recycling_space}. However, while for the orthogonal deflation preconditioner the effect of the ``low energy modes'' is to clearly reduce the average iteration number for a large portion of the recycle-restart plane, for the oblique deflation preconditioner it seems that there is little change in the lower-right half of the diagram, while on the upper-left half the iteration count is clearly increasing.
Consequently, we cannot draw clear conclusions regarding the optimal choice of SVD vectors.

\begin{figure}[tbp]
  \centering
  \includegraphics[width=0.49\columnwidth]{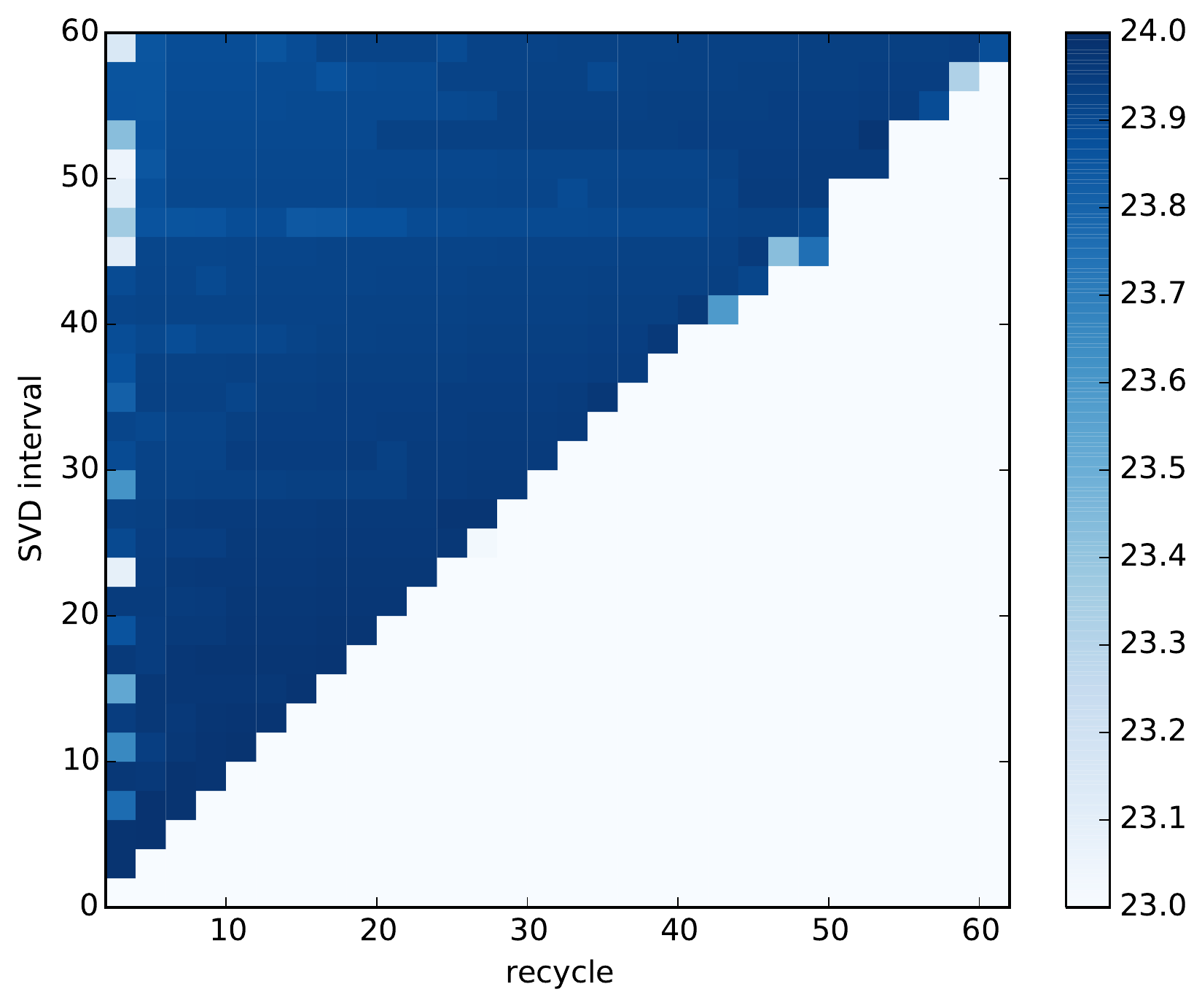}
  \includegraphics[width=0.49\columnwidth]{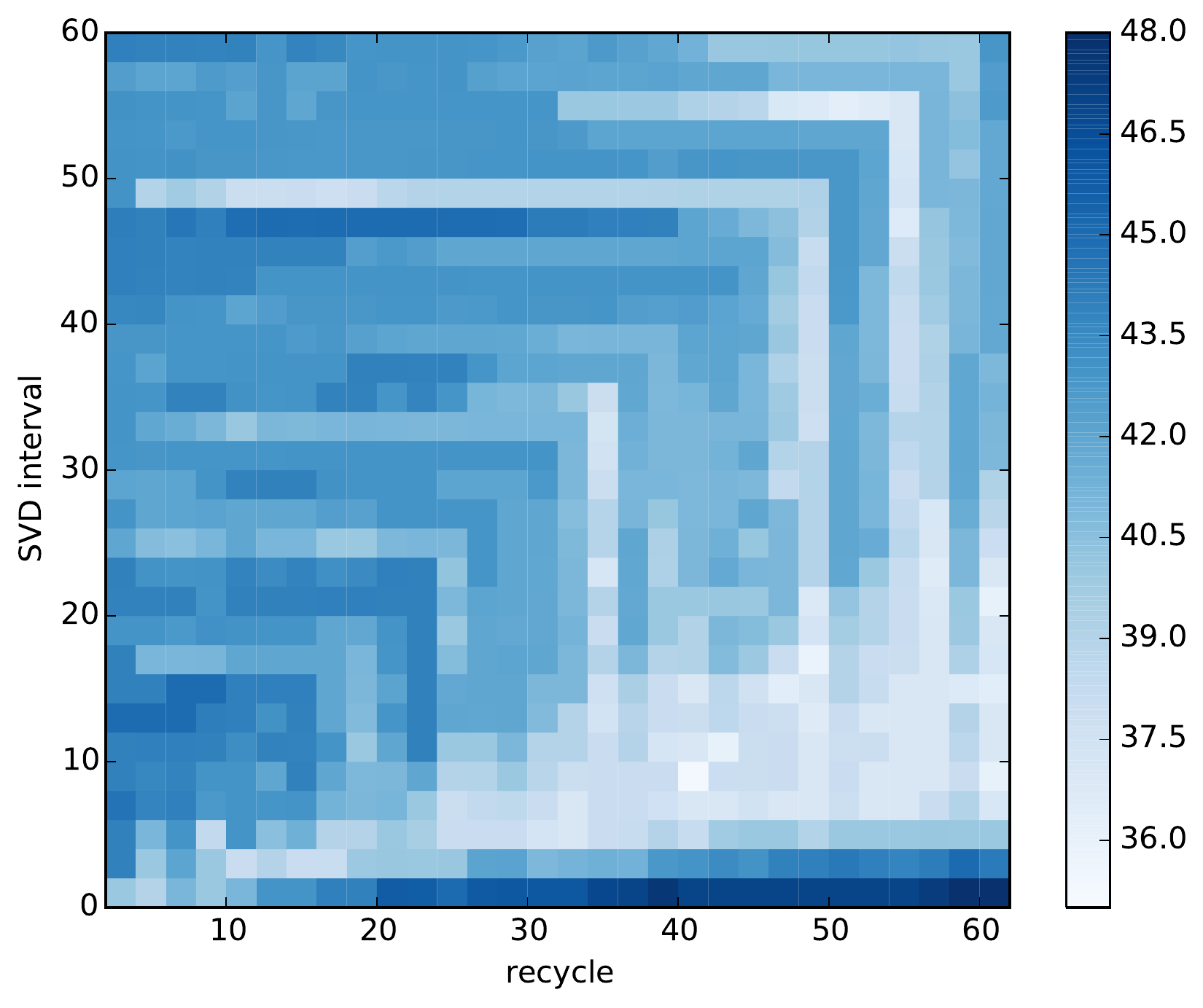}
  \caption{Average number of iterations to convergence for the \emph{augmented orthogonal} (left) and for the \emph{augmented oblique} (right) methods of table~\ref{tab:which_methods}. In this figure, the GMRES restart is fixed to 30, and the dimension of the recycling space and the interval between two SVD are varied.}
  \label{fig:map_aug_30}
\end{figure}

In all the numerical experiments carried out so far, the methods achieving the lowest iteration count have taken advantage of the SVD for the construction of the recycling spaces. However, the SVD is an expensive computation, whose cost more than balances the time saved by the fewer iterations needed to converge.
For this reason, we test the performance of the first two augmentation methods of Table~\ref{tab:which_methods}, except that in this case the SVD is not performed at each timestep, but only after a fixed number of timesteps, and the GMRES restart is fixed to 30.
The results of this set of simulations are available in Figure~\ref{fig:map_aug_30}. The maps in Figure~\ref{fig:map_aug_30}, although not directly comparable with the maps of the previous figures, show clearly that having a larger recycling space in this case can help to mitigate the effect of the fewer SVD selections carried out (and hence of the lower update frequency of the recycling spaces).

\subsection{Efficiency of SVD based preconditioners}
\label{sec:efficiency}
Despite the beneficial effects of the aforementioned recycling strategies on the iteration count and on the stability of the iterative solver, it is unlikely that such methods will be of interest if they fail in reducing the time required by the full computation.
For this reason, we give a back-of-the-envelope estimate on the opportunity to adopt recycling strategies.

\begin{figure}[htbp]
  \centering
  \includegraphics[width=0.8\textwidth]{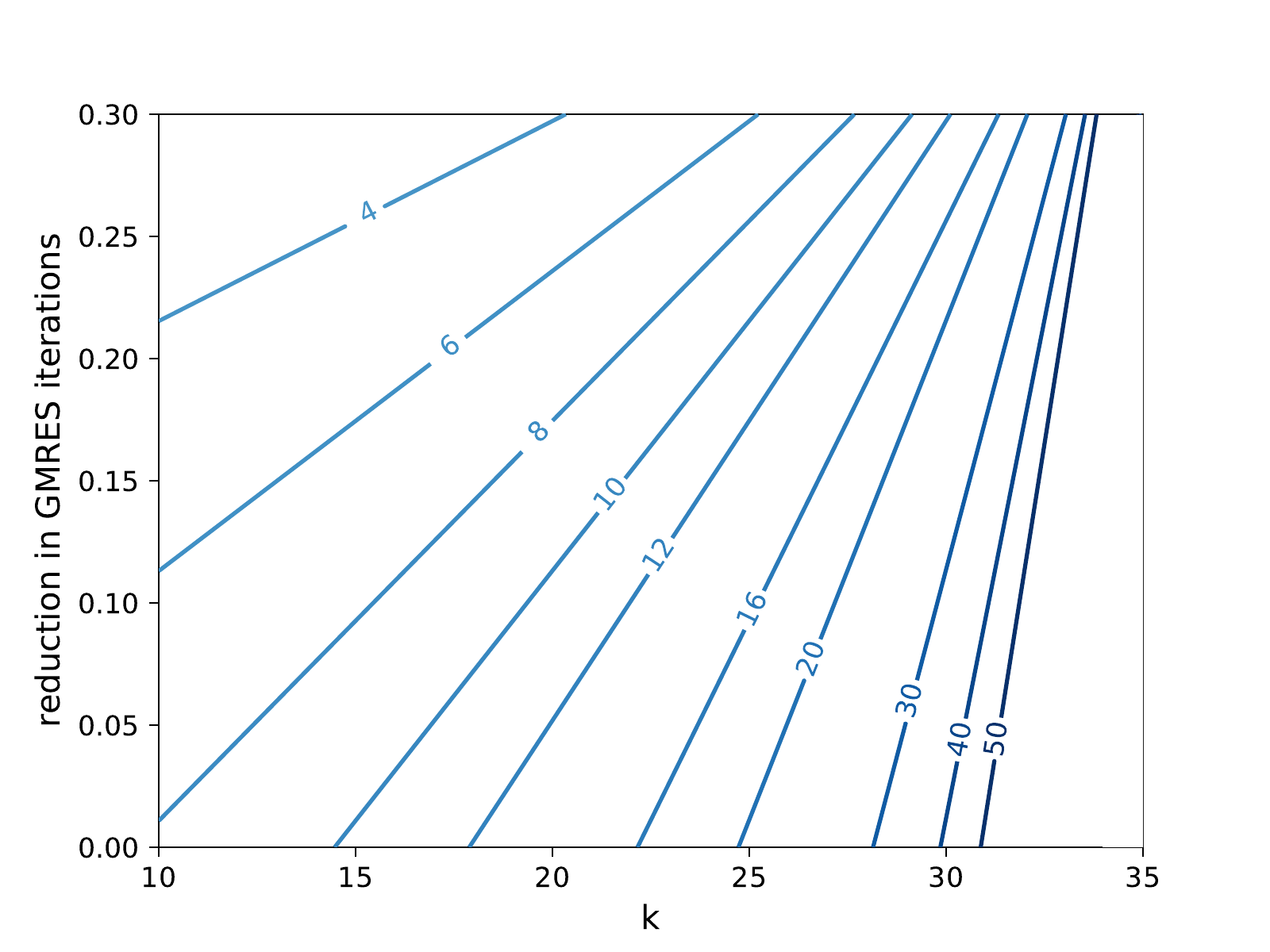}
  \caption{Minimum value of $\ell$ depending on the recycling space dimension $k$ and the ratio $\frac{r-\widetilde{r}}{r}$ in percentage.}
  \label{fig:plot_cost}
\end{figure}

Let $m$ be the dimension of the original linear system, $n$ the number of GMRES iterations before restarting, $b$ the bandwidth of the full matrix, $k$ the dimension of the recycling space. The operations required by the linear system solvers are then:
\begin{description}
  \item[dense GMRES] the Arnoldi orthogonalization algorithm requires $O(mn^2)$ operations, and the matrix-vector multiplications require $O(m^2n)$ operations. However, in the typical cases that lead to dense matrices (e.g. single domain Spectral Methods, Boundary Element Method), a fast matrix multiply is often available, reducing the cost of a matrix-vector product to $O(m\log m)$ or even to $O(m)$;
  \item[sparse GMRES] if the matrix is sparse, the cost of a matrix-vector multiplication reduces to $O(mb)$, and then the cost of each GMRES iteration reduces to $O(mb+mn)$;
  \item[truncated SVD] if $Z\in\mathbb{R}^{m\times s}$ is the matrix where the $s$ previous solutions are stored columnwise, then the computation of the first $k$ singular vectors from $Z$ requires $O(msk)+O(m^2k)$ operations~\cite{GolubVanLoan};
  \item[initial guess] the computation of an initial guess in the recycling space requires the solution of a linear system of dimension $k$, which can be achieved in $O(k^3)$ operations, and two matrix-vector multiplications, requiring $O(mk)$ operations;
  \item[augmented GMRES] in this case, at each iteration two rectangular matrix-vector products are required in addition to the usual matrix-vector product, for an overall cost of $O(mn^2)+O(m^2n)+O(mkn)$. For a sparse matrix, the cost is $O(mn^2)+O(mbn)+O(mkn)$.
\end{description}
Suppose now that the iteration count of the original system is $rn$, with $r$ number of GMRES restarts, and that for the recycling method it is $\widetilde{r}n$, with $\widetilde{r}<r$. The total cost $C$ and $C_r$ for each iteration of the system without recycling and for the system with recycling are respectively:
\begin{equation}
  \begin{aligned}
    C &= r\sum_{j=1}^n\left(mb+mj\right)\\
    C_4 &= \widetilde{r}\sum_{j=1}^n\left(mb+(j+k)m+(k+j)^3\right) +2mk+k^3+\frac{ksm+m^2k}{\ell},
  \end{aligned}
  \label{eq:total_costs}
\end{equation}
where we make the hypothesis that two consecutive SVD operations are spaced by $\ell$ timesteps. A reasonable way to proceed is to choose $n,s,k$, and $b$, make a guess on the ratio $r/\widetilde{r}$ and derive $\ell$ so that $C_r=C$. Then, any value of $\ell$ larger than this will lead to time savings\footnote{In modern cpu architectures, the direct proportionality between operation count and computational time is not strictly true, so this statement must be received with care.} if the ratio $r/\widetilde{r}$ does not change much with $\ell$.
For example, in Figure~\ref{fig:plot_cost} we show the minimum value of $\ell$ as a function of the recycling space dimension $k$ and the ratio $\frac{r-\widetilde{r}}{r}$ expressed in percentage. For the generation of Figure~\ref{fig:plot_cost}, we assumed $s=60$, $n=30$, $b=300$, $r=3$, $m=10^6$. Since the cost expressions~\eqref{eq:total_costs} are not exact (i.e. the coefficients are missing), this plot is valuable only as a trend indicator, and not as an exact tool for choice making. However, the trend in figure~\ref{fig:plot_cost} suggests that if the recycling method is not much effective in reducing the iteration count (lower part of the image), it is convenient to space the SVDs as much as possible, and that the same consideration holds as the recycling subspace dimension increases.
These considerations are highly dependent on the linear system dimension $m$: if the number of unknowns is very large, GMRES iterations will be very expensive, and even a small reduction in the iteration count can return significant time savings. Conversely, for small linear systems usually GMRES iterations are relatively inexpensive, and even a 30\% reduction in the iteration count can be easily overcome by the additional SVD expenditure.


To verify the time requirements of the \emph{augmentation orthogonal} method as opposed to the unrecycled case, we measure the time required to complete 1000 time steps of width $\Delta t=0.5$, with $\nu=10^{-2}$. In all the recycling runs, take an augmentation subspace of dimension 20, and we consider the following four choices for the generation of the recycling subspace:
\begin{description}
  \item[aug 20] the SVD is performed once every 20 time steps; the SVD is performed among the previous 20 solutions;
  \item[aug 40] the SVD is performed once every 40 time steps; the SVD is performed among the previous 40 solutions;
  \item[aug 60] the SVD is performed once every 60 time steps; the SVD is performed among the previous 60 solutions;
  \item[aug 80] the SVD is performed once every 80 time steps; the SVD is performed among the previous 80 solutions;
\end{description}
The results are reported in table~\ref{tab:comparison_svd} for the case of first order finite elements on a $64\times64$ grid. In terms of walltime, it seems that there is no advantage in choosing a recycling method over a traditional one. An obvious reason for this behaviour is the fact that first order finite elements lead to very sparse matrices, whose matrix-vector products are too inexpensive to produce any sensible time reduction, even with a 35\% reduction of the iteration count.

To test this hypothesis, we repeat the experiment with a $2\times2$ grid of 23\textsuperscript{rd} degree Legendre spectral elements. The related results are available in Table~\ref{tab:comparison_svd_sem}.

In all the cases of Tables~\ref{tab:comparison_svd} and~\ref{tab:comparison_svd_sem}, the auxiliary time required by the storage of previous solutions and by the SVD is less than 1\% of the solvers' time.

\begin{table}
\caption{Iterations and time required by the augmentation algorithm for first order finite elements. In the first row, the code is run without recycling preconditioner. 
In the remaining rows, the augmentation preconditioner is enabled using a recycling space of dimension 20 (for all rows), updated every 20, 40, 60 and 80 time steps respectively.}
\label{tab:comparison_svd}
\centering
\begin{tabular}{ccc}
\toprule
method & average iteration count & solver wall time (s) \\
\midrule
no aug & 124.5 & 222.7 \\
aug 20 & 78.36 & 249.0 \\
aug 40 & 78.00 & 243.6 \\
aug 60 & 77.48 & 268.9 \\
aug 80 & 78.20 & 271.8 \\
\bottomrule
\end{tabular}
\end{table}

\begin{table}
\caption{Iterations and time required by the augmentation algorithm for Legendre spectral elements of order 23. In the first row, the code is run without recycling preconditioner.
In the last four rows, the augmentation preconditioner is enabled using a recycling space of dimension 20 (for all rows), updated every 20, 40, 60 and 80 time steps respectively.}
\label{tab:comparison_svd_sem}
\centering
\begin{tabular}{ccc}
\toprule
method & average iteration count & solver wall time (s) \\
\midrule
no aug & 53.73 & 190.7 \\
aug 20 & 48.07 & 179.9 \\
aug 40 & 44.77 & 170.2 \\
aug 60 & 45.89 & 176.6 \\
aug 80 & 46.38 & 175.0 \\
\bottomrule
\end{tabular}
\end{table}

\section{Conclusions}
\label{sec:conclusions}

In this work we explored several strategies to accelerate the solution of slowly varying linear systems coming from the discretisation of time-dependent, non-linear partial differential equations with convection terms and diffusion. We developed a common framework for recycling methods that allowed to systematically construct deflative and augmentation methods, some of which are new.

The two families of recycling methods considered here were then compared for a choice of recycling subspaces coming from the SVD of previous solutions.
After introducing the new recycling subspaces, we studied their application in the construction of a good initial guess by means of Galerkin projection, as opposed to standard extrapolation.


Replacing extrapolation with projection introduces additional freedom in the choice of the projection subspace, and we studied the convergence properties of two different enriching subspaces: the first one was constructed by orthonormalizing only a few previous solution vectors, the second one by means of a Singular Value Decomposition of a larger set of previous solutions. This last approach is particularly convenient if one is interested anyway in constructing a Proper Orthogonal Decomposition approximation of some sort, and the POD basis is already available.


The recycling methods considered in this work were able to offer a reduction of the number of iterations in the GMRES solver in all the test cases. By weighting the overall cost of the algorithm, we conclude that SVD-accelerated methods may be an interesting strategy in all those cases where the matrices come from high order finite element methods, or spectral element methods, at least for scalar convection-diffusion problems, and moderate values of the P\'{e}clet number.

\bibliographystyle{abbrv}
\bibliography{deflation.bib}

\end{document}